\documentclass[10pt,letterpaper]{article}
\makeatletter
\usepackage{pdfsync}
\usepackage{cite}
\usepackage{amscd}
\usepackage{amsmath}
\usepackage{latexsym}
\usepackage{amsfonts}
\usepackage{amssymb}
\usepackage{amsthm}
\usepackage{graphicx}
\usepackage{indentfirst}
\usepackage{enumerate}
\usepackage{amstext}
\usepackage[dvipdfm,
            pdfstartview=FitH,
            CJKbookmarks=true,
            bookmarksnumbered=true,
            bookmarksopen=true,
            colorlinks=true, 
            pdfborder=001,   
            citecolor=magenta,
            linkcolor=blue,
            ]{hyperref}       

\newtheorem{thm}{\textbf Theorem}[section]
\newtheorem{lem}{\textbf Lemma}[section]
\newtheorem{rem}{\textbf Remark}[section]
\newtheorem{cor}{\textbf Corollary}[section]

\newtheorem*{thmA}{Theorem A}

\numberwithin{equation}{section}

\newcommand{\be}{\begin{eqnarray}}
\newcommand{\ee}{\end{eqnarray}}

\newcommand{\bes}{\begin{eqnarray*}}
\newcommand{\ees}{\end{eqnarray*}}

\begin{document}
\begin{titlepage}
\title{\bf On the propagation of regularity and decay of solutions to the Benjamin equation }

\author{Boling Guo$^{~a}$,  \quad Guoquan Qin$~^{b,*}$
\\[10pt]
\small {$^a $ Institute of Applied Physics and Computational Mathematics,
China Academy of Engineering Physics,}\\
\small {   Beijing,  100088,  P. R. China}\\[5pt]
\small {$^b $ Graduate School of China Academy of Engineering Physics,}\\
\small {  Beijing,  100088,  P. R. China}\\[5pt]
}
\footnotetext
{*~Corresponding author.

~~~~~E-mail addresses: gbl@iapcm.ac.cn(B. Guo), qinguoquan16@gscaep.ac.cn(G. Qin).}


\date{}
\end{titlepage}
\maketitle
\begin{abstract}
 In this paper, we investigate some special regularities
 and decay properties of  solutions to the initial value
 problem(IVP) of the Benjamin equation. The main result shows
 that: for initial datum $u_{0}\in H^{s}(\mathbb{R})$ with $s>3/4,$
 if the restriction of $u_{0}$ belongs to $H^{l}((x_{0}, \infty))$
  for some $l\in \mathbb{Z}^{+}$ and $x_{0}\in \mathbb{R},$
 then  the restriction of the corresponding solution $u(\cdot, t)$
 belongs to $H^{l}((\alpha, \infty))$ for any $\alpha\in \mathbb{R}$
 and any $t\in(0, T)$. Consequently, this type of regularity
 travels with infinite speed to its left as time evolves.


\end{abstract}

~~\noindent{ \textbf{Key words}: Benjamin equation; Propagation of regularity;
Decay}



\section{Introduction}
\setcounter{equation}{0}
In this paper, we are concerned with the IVP of the following Benjamin equation

\begin{eqnarray}\begin{cases}\label{benjamin}
u_{t}+\partial_{x}^{3}u-H\partial_{x}^{2}u+u\partial_{x}u=0,\quad x, t\in \mathbb{R},\\
u(x, 0)=u_{0}(x),
\end{cases}
\end{eqnarray}
where $H$ is the one-dimensional Hilbert transform

\begin{eqnarray*}\label{hilbert}
 Hf(x)&=&\frac{1}{\pi}\mbox{p.v.}\left(\frac{1}{x}*f\right)(x)\nonumber\\
 &=&\frac{1}{\pi}\lim_{\epsilon\rightarrow 0}\int_{|y|\geq \epsilon}
 \frac{f(x-y)}{y}\mbox{d}y\nonumber\\
 &=&(-i\mbox{sgn}(\xi)\hat{f}(\xi))^{\vee}(x)
 \end{eqnarray*}
 and $u=u(x, t)$ is a real valued function.

We will derive some special properties
including the  propagation of
regularity and decay of
solutions to equation
(\ref{benjamin}).

The integro-differential equation (\ref{benjamin}) models the unidirectional propagation of long
waves in a two-fluid system, where the lower fluid with greater density is infinitely deep
and the interface is subject to capillarity. It was derived by Benjamin \cite{be} to study gravity-capillary surface waves of solitary type on deep water. He also showed
that the solutions of the Benjamin equation (\ref{benjamin}) satisfy the following
conservation laws
\begin{eqnarray*}
 &&I_{1}(u)=\int_{-\infty}^{+\infty}u(x, t)\mbox{d}x,\\
 &&I_{2}(u)=\int_{-\infty}^{+\infty}u^{2}(x, t)\mbox{d}x,\\
 &&I_{3}(u)=\int_{-\infty}^{+\infty}
 [\frac{1}{2}(\partial_{x}u)^{2}(x, t)
 -\frac{1}{2}u(x, t)H\partial_{x}u(x, t)
 +\frac{1}{3}u^{3}(x, t)]\mbox{d}x.\\
\end{eqnarray*}
Notice that the conservation law for solutions of (\ref{benjamin})
\begin{equation*}
  I_{1}(u_{0})=\int_{-\infty}^{+\infty}u(x, t)\mbox{d}x
  =\int_{-\infty}^{+\infty}u_{0}(x)\mbox{d}x
\end{equation*}
guarantees that the property $\hat{u}(0)=0$ is preserved by the solution
flow.

Following  the definition of T. Kato  \cite{ka}
it is said that the  IVP (\ref{benjamin})
is locally well-posed (LWP) in the Banach space $X$
if given any datum $u_{0}\in X$ there exists $T>0$ and a unique solution
\begin{eqnarray}\label{class1}
u\in C([-T, T];X)\cap Y(T)
\end{eqnarray}
with $Y(T)$ be an auxillary function space.
Furthermore, the  solution map $u_{0}\mapsto u$ is continuous from
$X$ into the class (\ref{class1}).
This notion of LWP, which includes the ``persistent"
property, i.e., the solution
describes a continuous curve on $X$, implies that the
solution of (\ref{benjamin}) defines a  dynamic system on $X$.
If $T$ can be taken arbitrarily large, the IVP (\ref{benjamin}) is said to be
globally well-posed(GWP).


The problem of finding the minimal regularity property,
measured in the classical Sobolev space
\begin{equation*}
 H^{s}(\mathbb{R})=(1-\partial_{x}^{2})^{-s/2}L^{2}(\mathbb{R}),\quad s\in\mathbb{R},
\end{equation*}
required to guarantee that the  IVP
(\ref{benjamin}) is locally or globally well-posed in $H^{s}(\mathbb{R})$ has been
extensively studied.
We list some of the main results here.

Employing the Fourier restriction method introduced by Bourgain \cite{bo}
, Linares \cite{li} established
the LWP result for (\ref{benjamin}) in $H^{s}(\mathbb{R})$ with $s\geq 0,$
which combined with the conservation law $I_{2}$, leads to the GWP for (\ref{benjamin}) in $L^{2}.$
Guo and Huo \cite{gh2} obtained the LWP result  in $H^{s}(\mathbb{R})$ for
$s> -3/4.$
The best LWP results were established by Li and Wu \cite{lw}
 and Chen, Guo and Xiao \cite{cgx}. They also asserted the
 GWP for (\ref{benjamin})  in $H^{s}(\mathbb{R})$ for $s\geq -3/4.$
On the other hand, for the study of existence, stability and asymptotics of solitary wave solutions of
equation  (\ref{benjamin}), we can refer to \cite{be, abr, pa, cb, sl}.

The well-posedness problem has also been studied in the following
weighted Sobolev spaces  concerning with regularity and decay property
\begin{equation*}\label{weighted}
  Z_{s, r}=H^{s}(\mathbb{R})\cap L^{2}(|x|^{r}\mbox{d}x),\quad s, r\in \mathbb{R}
\end{equation*}
and
\begin{equation*}\label{weighted1}
 \dot{Z}_{s, r}=\{f\in Z_{s, r}: \hat{f}(0)=0\}.
\end{equation*}

In this respect
we can refer to, such as,  the articles  \cite{fp1, flp1, flp2}
for  the Benjamin-Ono and the dispersion generalized Benjamin-Ono equations,
the paper of Nahas and Ponce \cite{np} for the nonlinear
Schr\"{o}dinger equation, and so on.

For the Benjamin equation (\ref{benjamin}),
Urrea \cite{u} established the LWP  in weighted Sobolev spaces
$Z_{s, r}$ with $s\geq 1,$ $r\in[0, s/2]$  and $r<5/2$,
the GWP   in $Z_{s, r}$ with $s\geq 1,$ $r\in[0, s/2]$  and $3/2<r<5/2$,
 and the GWP in $\dot{Z}_{s, r}$ with $r\in[0, s/2]$  and $5/2\leq r<7/2$.
In particular, this implies the well-posedness of the  IVP (\ref{benjamin})
in the Schwartz space.
He also established a unique continuity property for  solutions of (\ref{benjamin}).
More precisely, he showed that if $u\in C([0, T]; Z_{7, 7/2^{-}})$ is a solution
of the IVP (\ref{benjamin}) and there exists three different times $t_{1}, t_{2}, t_{3}\in[0, T]$
such that $u(\cdot, t_{j})\in \dot{Z}_{7, 7/2}$ for $j=1, 2, 3,$ then $u(x, t)\equiv 0.$

Also, there are    works concerning with special regularities and decay properties
of some
dispersive models.

Isaza, Linares and Ponce \cite{ilp1} consider these problems for the
the k-generalized KdV equations
\begin{eqnarray}\begin{cases}\label{kdv}
u_{t}+\partial_{x}^{3}u+u^{k}\partial_{x}u=0,\quad x, t\in \mathbb{R},\\
u(x, 0)=u_{0}(x).
\end{cases}
\end{eqnarray}
They mainly established two results.

 The first one
describes the propagation of regularity in
the right hand side of the initial value for positive times.
It asserts that this regularity travels with infinite
speed to its left as time goes by.
Note that in \cite{ilp2}, they proved similar result for the following Benjamin-Ono
equation with negative dispersion
\begin{eqnarray}\begin{cases}\label{bo}
u_{t}-H\partial_{x}^{2}u+u\partial_{x}u=0,\quad x, t\in \mathbb{R},\\
u(x, 0)=u_{0}(x).
\end{cases}
\end{eqnarray}
 The difference between  \cite{ilp1} and  \cite{ilp2} lies in the
regularity of the initial data. For the
k-generalized KdV equations,  the initial value $u_{0}$ belongs to  $H^{3/4^{+}}(\mathbb{R})$,
while $u_{0}$ lies in  $H^{3/2}(\mathbb{R})$ for the Benjamin-Ono equation.

The second conclusion in \cite{ilp1}
is that
if the initial value $u_{0}\in H^{3/4^{+}}(\mathbb{R})$  of the k-generalized KdV equations has
polynomial decay in the positive real line,
then the corresponding solution possesses some persistence properties and
regularity effects for positive times.

Segata and   Smith  \cite{ss} extend the results of
\cite{ilp1} to the following fifth order dispersive equation
with $a_{1}, a_{2}, a_{3}$ be three constants
\begin{eqnarray}\begin{cases}\label{fifth}
u_{t}-\partial_{x}^{5}u+a_{1}u^{2}\partial_{x}u
+a_{2}\partial_{x}u\partial_{x}^{2}u
+a_{3}u\partial_{x}^{3}u=0,\quad x, t\in \mathbb{R},\\
u(x, 0)=u_{0}(x).
\end{cases}
\end{eqnarray}
However, the regularity of the initial data need to be $5/2^{+}$
for equation (\ref{fifth}).

Motivated by the above  works, the  objective of this paper
is to extend the results of  \cite{ilp1} to the IVP (\ref{benjamin}).

Before stating our results we
describe the following Theorem providing us with the space of solutions
where we shall be working on.

\begin{thmA}
Let $u_{0}\in H^{3/4^{+}}(\mathbb{R}).$ Then there exists a constant
$T=T(\|u_{0}\|_{H^{3/4^{+}}})$ and a unique local solution of
the  IVP (\ref{benjamin}) such that
\begin{eqnarray}\label{class}
 &&(i)\quad u\in C([-T, T]; H^{3/4^{+}}(\mathbb{R})),\nonumber\\
 &&(ii)\quad \partial_{x}u\in L^{4}([-T, T]; L^{\infty}(\mathbb{R})),\\
 &&(iii)\quad \sup_{x}\int_{-T}^{T}|J^{r}\partial_{x}u(x, t)|^{2}\mbox{d}t< \infty,
 \quad for \quad r\in[0,3/4^{+}],\nonumber\\
 &&(iv)\quad \int_{-\infty}^{\infty}\sup_{-T\leq t \leq T}|u(x, t)|^{2}\mbox{d}x< \infty,\nonumber
\end{eqnarray}
where $J=(1-\partial_{x}^{2})^{\frac{1}{2}}$ denotes the Bessel potential.
Moreover, the map data-solution, $u_{0}\mapsto u(x, t)$ is locally
continuous (smooth) from $H^{3/4^{+}}(\mathbb{R})$ into the class defined by
(\ref{class}).
\end{thmA}
\begin{rem}
The above well-posedness Theorem can  be
obtained by combining the properties of the unitary group
associated to the linear part of equation (\ref{benjamin}) and
the commutator estimate established by Kato and Ponce \cite{kp}.
For the method of its proof, we refer the reader to \cite{kpv} and \cite{la}, and
we omit the details here.
\end{rem}

We first describe the propagation of one-sided regularity
displayed by solutions to the IVP (\ref{benjamin}) provided by Theorem A.
\begin{thm}\label{regularity}
Assume $u_{0}\in H^{3/4^{+}}(\mathbb{R})$ and for some $l\in \mathbb{Z}^{+},$
 $l\geq 1$ and $x_{0}\in\mathbb{R}$ there holds
 \begin{eqnarray}
 \|\partial_{x}^{l}u_{0}\|_{L^{2}((x_{0}, \infty))}^{2}
 =\int_{x_{0}}^{\infty}|\partial_{x}^{l}u_{0}(x)|^{2}\mbox{d}x<\infty, \label{113}
 \end{eqnarray}
 then the solution of the  IVP (\ref{benjamin}) provided by
 Theorem A satisfies that for any $v>0$ and $\epsilon>0$
 \begin{eqnarray}\label{114}
\sup_{0\leq t\leq T}\int_{x_{0}+\epsilon-vt}^{\infty}(\partial_{x}^{j}u)^{2}(x, t)\mbox{d}x
\leq c_{0},
 \end{eqnarray}
 for $j=0, 1, 2, ..., l$ with $c_{0}=c_{0}(\|u_{0}\|_{H^{3/4^{+}}};  \|\partial_{x}^{l}u_{0}\|_{L^{2}((x_{0}, \infty))}; l; v; \epsilon; T)$.

 In particular, for all $t\in(0, T],$ the restriction of $u(\cdot, t)$ to any interval
 $(x_{1}, \infty)$ belongs to $H^{l}((x_{1}, \infty))$.

 Moreover, for any $v\geq0$ , $\epsilon>0 $ and $R>0$
 \begin{eqnarray}
 \int_{0}^{T}\int[D_{x}^{\frac{1}{2}}(\partial_{x}^{l}u(x, t)\eta(x+vt; \epsilon, b))]^{2}\mbox{d}x\mbox{d}t
 \leq c_{0},\label{115}\\
 \int_{0}^{T}\int_{x_{0}+\epsilon-vt}^{x_{0}+R-vt}
 (\partial_{x}^{l+1}u)^{2}(x, t)\mbox{d}x\mbox{d}t
\leq c_{1},\label{116}
 \end{eqnarray}
 where $c_{1}=c_{1}(l; \|u_{0}\|_{H^{3/4^{+}}};
\|\partial_{x}^{l}u_{0}\|_{L^{2}((x_{0}, \infty))}; v; \epsilon; T; R)$.

\end{thm}
\begin{rem}
The functions $\eta(x; \epsilon, b)$ mentioned in  Theorem \ref{regularity}
and $\eta_{j}(x; \epsilon, b)$ in Theorem \ref{persistence}
will be defined in section 2. In addition, without loss of generality, we shall
assume from now on $x_{0}=0$ in Theorem \ref{regularity}.
\end{rem}

The persistence of decay  and regularity effects established in
\cite{ilp1} can also be extended to the IVP (\ref{benjamin}).
In fact, we have

\begin{thm}\label{persistence}
Assume $u_{0}\in H^{3/4^{+}}(\mathbb{R})$ and for some $n\in \mathbb{Z}^{+},$
 $n\geq 1$  there holds
 \begin{eqnarray}
 \|x^{n/2}u_{0}\|_{L^{2}((0, \infty))}^{2}
 =\int_{0}^{\infty}|x^{n}||u_{0}(x)|^{2}\mbox{d}x<\infty, \label{117}
 \end{eqnarray}
 then the solution of the  IVP (\ref{benjamin}) provided by
 Theorem A satisfies that

 \begin{eqnarray}\label{118}
\sup_{0\leq t\leq T}\int_{0}^{\infty}|x^{n}||u(x, t)|^{2}\mbox{d}x
\leq c_{2}
 \end{eqnarray}
  with $c_{2}=c_{2}(\|u_{0}\|_{H^{3/4^{+}}}; \|x^{n/2}u_{0}\|_{L^{2}((0, \infty))}; T; n)$.

Furthermore, for any $v\geq0$ , $\epsilon, \delta>0 ,$ $m, j\in\mathbb{Z}^{+},$
 $m+j\leq n$ and  $m\geq 1,$

 \begin{eqnarray}
 &&\sup_{\delta\leq t\leq T}\int_{\epsilon-vt}^{\infty}(\partial_{x}^{m}u)^{2}(x, t)x_{+}^{j}\mbox{d}x
 +\int_{\delta}^{T}\int_{\epsilon-vt}^{\infty}
 (\partial_{x}^{m+1}u)^{2}(x, t)x_{+}^{j-1}\mbox{d}x\mbox{d}t\nonumber\\
 &&+\int_{\delta}^{T}\int[D_{x}^{\frac{1}{2}}(\partial_{x}^{m}u(x, t)\eta_{j}(x+vt; \epsilon, b))]^{2}\mbox{d}x\mbox{d}t
\leq c_{3},\label{119}
 \end{eqnarray}
 where  $c_{3}=c_{3}(l; \|u_{0}\|_{H^{3/4^{+}}}; \|x^{n/2}u_{0}\|_{L^{2}((x_{0}, \infty))}
; v; \epsilon; T; \delta; n),$ $x_{+}=\max\{x, 0\}.$
\end{thm}

Simple analysis of the proof
of Theorems \ref{regularity} and  \ref{persistence}
 yields their validity for the "defocusing"
 Benjamin equation
 \begin{eqnarray}\begin{cases}\label{dbenjamin}
u_{t}+\partial_{x}^{3}u-H\partial_{x}^{2}u-u\partial_{x}u=0,\quad x, t\in \mathbb{R},\\
u(x, 0)=u_{0}(x).
\end{cases}
\end{eqnarray}
Consequently, our results still hold
for $u(-x, -t)$ with $u(x, t)$ be the solution
of (\ref{benjamin}).
 Put another way,
for datum satisfying the assumption (\ref{113}) and
(\ref{117}) on the left hand side of the real line, respectively,
Theorems \ref{regularity} and  \ref{persistence} remain true backward in time.

On the other hand, equation (\ref{benjamin}) is time reversible.
In fact, let $v(x, t)=u(-x, -t)$ with $u(x, t)$ be
the solution of equation (\ref{benjamin}).
Using the relation $(Hv)(x, t)=-(Hu)(-x, -t)$,
one has
\begin{eqnarray}\begin{cases}\label{benjamin2}
v_{t}+\partial_{x}^{3}v-H\partial_{x}^{2}v+v\partial_{x}v=0,\quad x, t\in \mathbb{R},\\
v(x, 0)=u_{0}(-x).
\end{cases}
\end{eqnarray}
Theorems \ref{regularity} and  \ref{persistence}
combining with the above two points indicate
\begin{cor}\label{cor1}
Let $u\in C([-T, T]; H^{3/4^{+}}(\mathbb{R}))$ be a solution of
the equation (\ref{benjamin}) provided by Theorem A
such that
\begin{eqnarray*}
\partial_{x}^{m}u(\cdot, \hat{t})\notin L^{2}((a, \infty))\quad
for \quad some \quad \hat{t}\in (-T, T),\quad
 a\in \mathbb{R}\quad and \quad m\in \mathbb{Z}^{+}.
\end{eqnarray*}
Then for any $t\in [-T, \hat{t})$ and any $\beta\in \mathbb{R}$
\begin{equation*}
\partial_{x}^{m}u(\cdot, t)\notin L^{2}((\beta, \infty)),
\quad and \quad x^{m/2}u(\cdot, t)\notin L^{2}((0, \infty)).
\end{equation*}
\end{cor}
Next, Theorems \ref{regularity} and  \ref{persistence} yield that
the singularity of the solution corresponding to an
appropriate class of initial data propagates with infinite
speed  to the left as time goes by.
Also, since equation (\ref{benjamin}) is time
reversible, the solution cannot have had some
regularity in the past.
More precisely, we have
\begin{cor}\label{cor2}
Let $u\in C([-T, T]; H^{3/4^{+}}(\mathbb{R}))$ be a solution of
the equation (\ref{benjamin}) provided by Theorem A.
Suppose there exists  $n, m\in \mathbb{Z}^{+}$ with
$m\leq n$ such that for some $a, b\in \mathbb{R}$
with $a< b$
\begin{eqnarray}\label{cor119}
\int_{b}^{\infty}|\partial_{x}^{n}u_{0}(x)|^{2}\mbox{d}x<\infty \quad
but \quad
\partial_{x}^{m}u_{0}\notin L^{2}((a, \infty)).
\end{eqnarray}

Then for any $t\in (0, T)$ and any $v, \epsilon> 0$
\begin{equation*}
\int_{b+\epsilon-vt}^{\infty}|\partial_{x}^{n}u(x, t)|^{2}\mbox{d}x<\infty
\end{equation*}
and for any $t\in (-T, 0)$ and any $\alpha\in \mathbb{R}$
\begin{equation*}
\int_{\alpha}^{\infty}|\partial_{x}^{m}u(x, t)|^{2}\mbox{d}x=\infty.
\end{equation*}
\end{cor}

We now discuss some of the ingredients in the proof of
Theorems \ref{regularity} and \ref{persistence}.

The first one is concerned with the proof of
Theorem \ref{regularity}. As in \cite{ilp1}, we mainly use induction.
To treat the Benjamin-Ono term $-H\partial_{x}^{2}u,$
we follow the idea in \cite{ilp2}, where the commutator
estimate for the Hilbert transform (\ref{commutator}) plays
a vital role.
In spite of this, there is a little difference between
\cite{ilp2}  and this paper
when handling the following two terms (see (\ref{2A422}))
\begin{eqnarray}
\int_{0}^{T}\int(\partial_{x}^{2}u)^{2} (\eta^{'})^{2}\mbox{d}x\mbox{d}t
+\int_{0}^{T}\int(\partial_{x}^{2}u \eta)^{2}\mbox{d}x\mbox{d}t.\label{q1}
\end{eqnarray}
In \cite{ilp2} for the Benjamin-Ono equation (\ref{bo}), these two terms
can be controlled by sufficient local smoothing effect.
More precisely, the condition (1.6)(iii) in \cite{ilp2} reads
\begin{eqnarray*}
\int_{-T}^{T}\int_{-R}^{R}
(|\partial_{x}D_{x}u|^{2}+|\partial_{x}^{2}u|^{2})\mbox{d}x\mbox{d}t\leq c_{0},
\end{eqnarray*}
where $R$ is arbitrary and finite.
This combined with the boundedness of $\eta^{'}$ and $\eta$ on the support
of $\eta$ immediately yields the finiteness of (\ref{q1}).
However, (\ref{class})(iii) in this paper provides us
at most $7/4^{+}$ order local smoothing effect, which is not enough
to bound (\ref{q1}).
Fortunately, in the first step(the case l=1 in the proof of
Theorem \ref{regularity}) in our
induction process, the KdV term provides us with the finiteness of (see (\ref{conclusion1}))
\begin{eqnarray*}
  \int_{0}^{T}\int (\partial_{x}^{2}u)^{2}(x, t)\chi^{'}(x+vt; \epsilon, b)\mbox{d}x\mbox{d}t.\label{q2}
\end{eqnarray*}
This permits us to use the properties of $\eta^{'}$ and $\eta$,
i.e. (\ref{prp3}) and (\ref{prp4}), to control (\ref{q1}).

The second one relates to the proof of
Theorem \ref{persistence}.
The difficulty still comes from the Benjamin-Ono term.
For the term $A_{422}$ in (\ref{4A42}), note that because of the factor $x^{n}$
in the definition of $\eta_{n}^{'}$ and $\eta_{n}$,
the support of $\eta_{n}$ is not $[\epsilon, b]$ at all for the general case.
As a consequence,   $\eta_{n}$
and $\eta_{n}^{'}$ may be unbounded.
However, we notice that (\ref{prp10}) and (\ref{prp11})
provide us with a relation between $\chi_{n} $ and $\chi_{n-1} $,
therefore, we could use induction to treat this term.
(\ref{prp10}) and (\ref{prp11}) are also used to bound the term in
(\ref{6A422}).


The rest of this paper is organized as follows: in section 2
we construct  our cut-off functions and state a lemma to be used in
the proof of Theorem \ref{regularity} and  Theorem \ref{persistence}.
The proof of Theorem \ref{regularity} and  Theorem \ref{persistence}
will be given in section 3 and section 4, respectively.

\section{Preliminaries}
Let us  first  construct our cutoff functions,
the construction of this family of cutoff functions
is motivated by Segata and Smith \cite{ss}.


Let $p$ be large enough and let $\rho(x)$ be defined as follows
\begin{eqnarray*}
\rho(x)=a\int_{0}^{x}y^{p}(1-y)^{p}\mbox{d}y
\end{eqnarray*}
with the constant $a=a(p)$ be chosen to  satisfy  $\rho(1)=1.$

\begin{rem}
According to Lemma \ref{lem1} below, when come across the
$L^{p}$ norm of the commutator related to the Hilbert
transform, we want to put all derivatives to the smooth function
$\psi,$ this is the reason for $p$ in the definition
of $\rho(x)$ being large enough.
 \end{rem}

With the above definition, we have
\begin{eqnarray*}
\rho(0)=0,\quad \rho(1)=1,\\
\rho^{'}(0)=\rho^{''}(0)=\cdot\cdot\cdot=\rho^{(p)}(0)=0,\\
\rho^{'}(1)=\rho^{''}(1)=\cdot\cdot\cdot=\rho^{(p)}(1)=0
\end{eqnarray*}
with $0<\rho, \rho^{'}$ for $0<x<1. $

Next, for parameters $\epsilon, b>0,$ define $\chi\in C^{p}(\mathbb{R})$ by
\begin{eqnarray*}
\chi(x; \epsilon, b)=
\begin{cases}
0,\quad x\leq \epsilon\\
\rho((x-\epsilon)/b),\quad \epsilon<x<b+\epsilon\\
1,\quad b+\epsilon\leq x.
\end{cases}
\end{eqnarray*}
In addition, we define $\chi_{n}=x^{n}\chi\in C^{p}(\mathbb{R}).$

By their definitions, $\chi $ and $\chi_{n}$ are both  positive for $x\in (\epsilon, \infty)$.

  Computing as section 2 in \cite{ss}, we can derive  the following properties concerning
  $\chi $ and $\chi_{n}$:
\begin{eqnarray}
&&(1)\quad \chi(x; \epsilon/10, \epsilon/2 )=1\quad on \quad\mbox{supp}\chi(x; \epsilon, b )=[\epsilon, \infty);\label{id}\\
&&(2)\quad |\chi(x, \epsilon, b)|\leq \chi_{1}^{'}(x, \epsilon, b);\label{prp2}\\
&&(3)\quad \left|\frac{[\chi^{''}(x; \epsilon, b)]^{2}}{\chi^{'}(x; \epsilon, b)}\right|
\leq c(\epsilon, b)\chi^{'}(x; \epsilon/3, b+\epsilon) \quad on \quad support \quad of \quad\chi^{'};\label{prp3}\\
&&(4)\quad |\chi^{(j)}(x; \epsilon, b)|
\leq c(j, \epsilon, b)\chi^{'}(x; \epsilon/3, b+\epsilon)\quad on \quad[\epsilon, b+\epsilon]\quad
 for \quad j=1, 2,..., p;\label{prp4}\\
&&(5)\quad |\chi_{n-l}^{''}(x, \epsilon, b)|\leq c(n, l)\chi_{n-l-2}(x, \epsilon, b)\nonumber\\
&&\quad \quad+c(b, v, \epsilon, T)\chi^{'}(x, \epsilon/3, b+\epsilon)
\quad for \quad l\leq n-2;\label{prp5}\\
&&(6)\quad |\chi_{n-l}^{'''}(x, \epsilon, b)|
\leq c(n, l)\chi_{n-l-3}(x, \epsilon, b)\nonumber\\
&&\quad \quad+c(n, l, b)\chi(x, \epsilon/10,\epsilon/2)\quad
for \quad l\leq n-3;\label{prp6}\\
&&(7)\quad n\chi_{n-1}(x, \epsilon, b)\leq\chi_{n}^{'}(x, \epsilon, b) ;\label{prp7}\\
&&(8)\quad |\chi_{n}^{(j)}(x; \epsilon, b)|
\leq c(j, n, b)[1+\chi_{n}(x; \epsilon, b)]\quad
for \quad j=1, 2, 3, ..., p;\label{prp8}\\
&&(9)\quad |\chi_{n}^{(j)}(x; \epsilon, b)|
\leq c(\epsilon, n, b)\chi_{n-1}(x; \epsilon/3, b+\epsilon)\quad
for \quad j=1, 2, 3, ..., p;\label{prp10}\\
&&(10)\left|\frac{[\chi_{n}^{''}(x; \epsilon, b)]^{2}}{\chi_{n}^{'}(x; \epsilon, b)}\right|
\leq c(\epsilon, b, n)\chi_{n-1}(x; \epsilon/3, b+\epsilon)\quad on \quad support \quad of \quad\chi_{n}^{'}.\label{prp11}
\end{eqnarray}

Moreover, we define
\begin{eqnarray}\label{eta}
  \eta(x; \epsilon, b)=\sqrt{\chi^{'}(x; \epsilon, b )},\nonumber\\
  \eta_{n}(x; \epsilon, b)=\sqrt{\chi_{n}^{'}(x; \epsilon, b )}.
\end{eqnarray}
Then, reasoning as section 2 in \cite{ilp2}, we derive that $\eta(x; \epsilon, b)$ and $\eta_{n}(x; \epsilon, b)$
are both in $C^{p}(\mathbb{R}).$

The following commutator estimate is an extension of the Calder\'{o}n
theorem \cite{ca}, it was proved by Dawson, McGahagan and Ponce \cite{dmp}.
\begin{lem}\label{lem1}
For any $p\in (1, \infty)$ and $l, m\in\mathbb{Z}^{+}\cup \{0\},$ $l+m\geq 1$
there exists a constant $C=C(p, l, m)>0$ such that
\begin{equation}\label{commutator}
  \|\partial_{x}^{l}[H; \psi]\partial_{x}^{m}f\|_{L^{p}}
  \leq C\|\partial_{x}^{l+m}\psi\|_{L^{\infty}}\|f\|_{L^{p}}
\end{equation}
with $H$ be the Hilbert transform.
\end{lem}

\section{Proof of Theorem \ref{regularity} }

To prove Theorem \ref{regularity}, we follow the idea in \cite{ilp1} and use an induction argument.
To  illuminate our method, we first prove (\ref{114}) for $l=1$ and $l=2$.

Let us first prove the case $l=1.$

Formally, applying  $\partial_{x}$ to equation (\ref{benjamin})
and multiplying the result by $\partial_{x}u(x, t)\chi(x+vt; \epsilon, b)$,
after some integration by parts, one deduces
\begin{eqnarray}\label{first}
&&\frac{1}{2}\frac{\mbox{d}}{\mbox{d}t}\int(\partial_{x}u)^{2}(x, t)\chi(x+vt)\mbox{d}x
\underbrace{-\frac{1}{2}v\int(\partial_{x}u)^{2}(x, t)\chi^{'}(x+vt)\mbox{d}x}_{
\text{$A_{1}$}}\nonumber\\
&&+\frac{3}{2}v\int(\partial_{x}^{2}u)^{2}(x, t)\chi^{'}(x+vt)\mbox{d}x
\underbrace{-\frac{1}{2}\int(\partial_{x}u)^{2}(x, t)\chi^{'''}(x+vt)\mbox{d}x}_{
\text{$A_{2}$}}\nonumber\\
&&\underbrace{+\int\partial_{x}(u\partial_{x}u)\partial_{x}u(x, t)\chi(x+vt)\mbox{d}x}_{
\text{$A_{3}$}}
\underbrace{-\int H\partial_{x}^{3}u\partial_{x}u(x, t)\chi(x+vt)\mbox{d}x}_{
\text{$A_{4}$}}=0,
\end{eqnarray}
where in $\chi$ we omit the parameters $\epsilon$ and $b$.

We  estimate the integrals in (\ref{first})  term by term.

Using (\ref{class})(iii) with $r=0$ and the support property of
$\chi^{'}(x)$ , it holds that
\begin{eqnarray*}
\int_{0}^{T}|A_{1}(t)|\mbox{d}t
\leq \int_{0}^{T}\int(\partial_{x}u)^{2}(x, t)\chi^{'}(x+vt)\mbox{d}x\mbox{d}t\leq c_{0}
\end{eqnarray*}
and similarly
\begin{eqnarray*}
\int_{0}^{T}|A_{2}(t)|\mbox{d}t
\leq c_{0}.
\end{eqnarray*}
For the term $A_{3}$, direct computation yields
\begin{eqnarray*}
A_{3}(t)
&=&\int(\partial_{x}u)^{3}\chi(x+vt)\mbox{d}x
+\int u \partial_{x}^{2}u\partial_{x}u\chi(x+vt)\mbox{d}x\nonumber\\
&=&\frac{1}{2}\int(\partial_{x}u)^{3}\chi(x+vt)\mbox{d}x
-\frac{1}{2}\int u \partial_{x}u\partial_{x}u\chi^{'}(x+vt)\mbox{d}x\nonumber
\end{eqnarray*}
\begin{eqnarray*}
\quad\quad\quad\quad\quad\quad&\leq& \|\partial_{x}u\|_{L^{\infty}}\int(\partial_{x}u)^{2}\chi(x+vt)\mbox{d}x
+\|u\|_{L^{\infty}}\int(\partial_{x}u)^{2}\chi^{'}(x+vt)\mbox{d}x\nonumber\\
&=&A_{31}+A_{32}.
\end{eqnarray*}
By Sobolev embedding theorem, one obtains
\begin{eqnarray*}
\int_{0}^{T}|A_{32}(t)|\mbox{d}t
\leq\sup_{[0, T]}\|u\|_{H^{3/4^{+}}}
\int(\partial_{x}u)^{2}\chi^{'}(x+vt)\mbox{d}x\mbox{d}t.
\end{eqnarray*}
The term $A_{31}$ will be controlled by using
(\ref{class})(ii) and the Gronwall inequality.

Finally, to estimate $A_{4}$, we follow the idea described in \cite{ilp2}.

Integration by parts yields
\begin{eqnarray*}
A_{4}
&=&-\int H\partial_{x}^{3}u\partial_{x}u\chi(x+vt)\mbox{d}x\nonumber\\
&=&\int H\partial_{x}^{2}u\partial_{x}^{2}u\chi(x+vt)\mbox{d}x
+\int H\partial_{x}^{2}u\partial_{x}u\chi^{'}(x+vt)\mbox{d}x\nonumber\\
&=&A_{41}+A_{42}.
\end{eqnarray*}
Since the Hilbert transform is skew symmetric,
we have
\begin{eqnarray*}
A_{41}
&=&\int H\partial_{x}^{2}u\partial_{x}^{2}u\chi(x+vt)\mbox{d}x\nonumber\\
&=&-\int \partial_{x}^{2}u H(\partial_{x}^{2}u\chi(x+vt))\mbox{d}x\nonumber\\
&=&-\int \partial_{x}^{2}u H\partial_{x}^{2}u\chi(x+vt)\mbox{d}x
-\int \partial_{x}^{2}u [H; \chi]\partial_{x}^{2}u\mbox{d}x\nonumber\\
&=&-A_{41}-\int \partial_{x}^{2}u [H; \chi]\partial_{x}^{2}u\mbox{d}x.
\end{eqnarray*}
Therefore,  (\ref{commutator}) leads to
\begin{eqnarray*}
A_{41}
&=&-\frac{1}{2}\int \partial_{x}^{2}u [H; \chi]\partial_{x}^{2}u\mbox{d}x\nonumber\\
&=&-\frac{1}{2}\int u \partial_{x}^{2}[H; \chi]\partial_{x}^{2}u\mbox{d}x\nonumber\\
&\leq& c\|u\|_{L^{2}}\|\partial_{x}^{2}[H; \chi]\partial_{x}^{2}u\|_{L^{2}}\nonumber\\
&\leq& c\|u\|_{L^{2}}^{2}=c\|u_{0}\|_{L^{2}}^{2}.
\end{eqnarray*}
Concerning  the term $A_{42}$, let us recall the definition of $\eta(x; \epsilon, b)$
in (\ref{eta}), we can write $A_{42}$ as
\begin{eqnarray*}
A_{42}
&=&\int H\partial_{x}^{2}u\partial_{x}u\chi^{'}(x+vt)\mbox{d}x\nonumber\\
&=&\int H\partial_{x}^{2}u \eta \partial_{x}u\eta\mbox{d}x\nonumber\\
&=&\int  H(\partial_{x}^{2}u \eta) \partial_{x}u\eta\mbox{d}x
-\int [H; \eta]\partial_{x}^{2}u\partial_{x}u\eta\mbox{d}x\nonumber\\
&=&\int  H\partial_{x}(\partial_{x}u \eta) \partial_{x}u\eta\mbox{d}x
-\int  H(\partial_{x}u \eta^{'}) \partial_{x}u\eta\mbox{d}x
-\int [H; \eta]\partial_{x}^{2}u\partial_{x}u\eta\mbox{d}x\nonumber\\
&=&A_{421}+A_{422}+A_{423}.
\end{eqnarray*}
Plancherel's identity yields
\begin{eqnarray}\label{1A421}
A_{421}=\int  H\partial_{x}(\partial_{x}u \eta) \partial_{x}u\eta\mbox{d}x
=\int[D_{x}^{\frac{1}{2}}(\partial_{x}u \eta)]^{2}\mbox{d}x,
\end{eqnarray}
which is positive and will stay  at the left hand side of (\ref{first}).

 The  boundedness of the Hilbert transform in $L^{2}$
 and the Young inequality produce
\begin{eqnarray}\label{A422}
\int_{0}^{T}|A_{422}|\mbox{d}t
&=&\int_{0}^{T}\left|\int  H(\partial_{x}u \eta^{'}) \partial_{x}u\eta\mbox{d}x\right|\mbox{d}t\nonumber\\
&\leq& c\int_{0}^{T}\int(\partial_{x}u)^{2} (\eta^{'})^{2}\mbox{d}x\mbox{d}t
+c\int_{0}^{T}\int(\partial_{x}u)^{2}\eta^{2}\mbox{d}x\mbox{d}t.
\end{eqnarray}
Employing  the boundedness of $\eta$ and $\eta^{'}$ on the support of $\eta$
and using (\ref{class})(iii) with $r=0$, we obtain
\begin{eqnarray*}
\int_{0}^{T}|A_{422}|\mbox{d}t\leq c_{0}.
\end{eqnarray*}
Invoking the commutator estimate (\ref{commutator}), we derive
\begin{eqnarray*}
|A_{423}|
&=&\left|\int [H; \eta]\partial_{x}^{2}u\partial_{x}u\eta\mbox{d}x\right|\nonumber\\
&\leq& \|[H; \eta]\partial_{x}^{2}u\|_{L^{2}}\|\partial_{x}u\eta\|_{L^{2}}\nonumber\\
&\leq& c\|u\|_{L^{2}}^{2}+c\|\partial_{x}u\eta\|_{L^{2}}^{2}\nonumber\\
&\leq& c\|u_{0}\|_{L^{2}}^{2}+c\|\partial_{x}u\eta\|_{L^{2}}^{2}.
\end{eqnarray*}
After integration in time, the term $\|\partial_{x}u\eta\|_{L^{2}}^{2}$
can be controlled as that in (\ref{A422}).

Substituting  the above information in $(\ref{first})$, using the
Gronwall inequality and (\ref{class})(ii), one obtains

\begin{eqnarray}\label{conclusion1}
&&\sup_{t\in[0, T]}\int (\partial_{x}u)^{2}(x, t)\chi(x+vt; \epsilon, b)\mbox{d}x
+\int_{0}^{T}\int (\partial_{x}^{2}u)^{2}(x, t)\chi^{'}(x+vt; \epsilon, b)\mbox{d}x\mbox{d}t\nonumber\\
&&+\int_{0}^{T}\int [D_{x}^{\frac{1}{2}}(\partial_{x}u(x, t)
 \eta(x+vt; \epsilon, b))]^{2}\mbox{d}x\mbox{d}t
\leq c_{0}.
\end{eqnarray}
This completes the proof of the case $l=1.$

Next, we prove (\ref{114}) for the case $l=2.$

Applying $\partial_{x}^{2}$ to equation (\ref{benjamin}), then  multiplying
$\partial_{x}^{2}u(x, t)\chi(x+vt; \epsilon, b)$ and integrating, we find
\begin{eqnarray}\label{second}
&&\frac{1}{2}\frac{\mbox{d}}{\mbox{d}t}\int(\partial_{x}^{2}u)^{2}(x, t)\chi(x+vt)\mbox{d}x
\underbrace{-\frac{1}{2}v\int(\partial_{x}^{2}u)^{2}(x, t)\chi^{'}(x+vt)\mbox{d}x}
_{\text{$A_{1}$}}\nonumber\\
&&+\frac{3}{2}v\int(\partial_{x}^{3}u)^{2}(x, t)\chi^{'}(x+vt)\mbox{d}x
\underbrace{-\frac{1}{2}\int(\partial_{x}^{2}u)^{2}(x, t)\chi^{'''}(x+vt)\mbox{d}x}
_{\text{$A_{2}$}}\nonumber\\
&&\underbrace{+\int\partial_{x}^{2}(u\partial_{x}u)\partial_{x}^{2}u(x, t)\chi(x+vt)\mbox{d}x}
_{\text{$A_{3}$}}
\underbrace{-\int H\partial_{x}^{4}u\partial_{x}^{2}u(x, t)\chi(x+vt)\mbox{d}x}
_{\text{$A_{4}$}}=0.
\end{eqnarray}

 Invoking (\ref{conclusion1}), one has
 \begin{eqnarray*}
\int_{0}^{T}|A_{1}(t)|\mbox{d}t
\leq |v|\int_{0}^{T}\int(\partial_{x}^{2}u)^{2}(x, t)\chi^{'}(x+vt)\mbox{d}x
\leq c_{0}.
\end{eqnarray*}
Employing (\ref{prp4}) with $j=3$ and using (\ref{conclusion1})
with $(\epsilon/3, b+\epsilon)$ instead of $(\epsilon, b)$, it holds that
\begin{eqnarray}\label{2A2}
\int_{0}^{T}|A_{2}(t)|\mbox{d}t
&\leq&\int_{0}^{T}\int(\partial_{x}^{2}u)^{2}(x, t)|\chi^{'''}(x+vt; \epsilon, b)|\mbox{d}x\mbox{d}t\nonumber\\
&\leq& \int_{0}^{T}\int(\partial_{x}^{2}u)^{2}(x, t)|\chi^{'}(x+vt; \epsilon/3, b+\epsilon)|\mbox{d}x\mbox{d}t
\leq c_{0}.
\end{eqnarray}

Integration by parts yields
\begin{eqnarray*}
A_{3}(t)
&=&3\int\partial_{x}u(\partial_{x}^{2}u)^{2}\chi(x+vt)\mbox{d}x
+\int u \partial_{x}^{3}u\partial_{x}^{2}u\chi(x+vt)\mbox{d}x\nonumber\\
&=&\frac{5}{2}\int\partial_{x}u(\partial_{x}^{2}u)^{2}\chi(x+vt)\mbox{d}x
-\frac{1}{2}\int u (\partial_{x}^{2}u)^{2}\chi^{'}(x+vt)\mbox{d}x\nonumber\\
&\leq& \|\partial_{x}u\|_{L^{\infty}}\int(\partial_{x}^{2}u)^{2}\chi(x+vt)\mbox{d}x
+\|u\|_{L^{\infty}}\int(\partial_{x}^{2}u)^{2}\chi^{'}(x+vt)\mbox{d}x\nonumber\\
&=&A_{31}+A_{32}.
\end{eqnarray*}
Again, using Sobolev embedding theorem and (\ref{conclusion1}), one obtains
\begin{eqnarray*}
\int_{0}^{T}|A_{32}(t)|\mbox{d}t
<\sup_{[0, T]}\|u\|_{H^{3/4^{+}}}
\int_{0}^{T}\int(\partial_{x}^{2}u)^{2}\chi^{'}(x+vt)\mbox{d}x\mbox{d}t\leq c_{0}.
\end{eqnarray*}
The term $A_{31}$ will be controlled by using
(\ref{class})(ii) and the Gronwall inequality.

We now estimate $A_{4}$.

Integration by parts leads to
\begin{eqnarray*}
A_{4}
&=&-\int H\partial_{x}^{4}u\partial_{x}^{2}u\chi(x+vt)\mbox{d}x\nonumber\\
&=&\int H\partial_{x}^{3}u\partial_{x}^{3}u\chi(x+vt)\mbox{d}x
+\int H\partial_{x}^{3}u\partial_{x}^{2}u\chi^{'}(x+vt)\mbox{d}x\nonumber\\
&=&A_{41}+A_{42}.
\end{eqnarray*}
Invoking again the fact that the Hilbert transform is skew symmetric,
we have
\begin{eqnarray}\label{2A41}
A_{41}
&=&\int H\partial_{x}^{3}u\partial_{x}^{3}u\chi(x+vt)\mbox{d}x\nonumber\\
&=&-\int \partial_{x}^{3}u H(\partial_{x}^{3}u\chi(x+vt))\mbox{d}x\nonumber\\
&=&-\int \partial_{x}^{3}u H\partial_{x}^{3}u\chi(x+vt)\mbox{d}x
-\int \partial_{x}^{3}u [H; \chi]\partial_{x}^{3}u\mbox{d}x\nonumber\\
&=&-A_{41}-\int \partial_{x}^{3}u [H; \chi]\partial_{x}^{3}u\mbox{d}x.
\end{eqnarray}
Consequently, (\ref{commutator}) produces
\begin{eqnarray*}
A_{41}
&=&-\frac{1}{2}\int \partial_{x}^{3}u [H; \chi]\partial_{x}^{3}u\mbox{d}x\nonumber\\
&=&\frac{1}{2}\int u \partial_{x}^{3}[H; \chi]\partial_{x}^{3}u\mbox{d}x\nonumber\\
&\leq& c\|u\|_{L^{2}}\|\partial_{x}^{3}[H; \chi]\partial_{x}^{3}u\|_{L^{2}}\nonumber\\
&\leq& c\|u\|_{L^{2}}^{2}=c\|u_{0}\|_{L^{2}}^{2}.
\end{eqnarray*}

Applying (\ref{eta}) yields
\begin{eqnarray*}
A_{42}
&=&\int H\partial_{x}^{3}u\partial_{x}^{2}u\chi^{'}(x+vt)\mbox{d}x\nonumber\\
&=&\int H\partial_{x}^{3}u \eta \partial_{x}^{2}u\eta\mbox{d}x\nonumber\\
&=&\int  H(\partial_{x}^{3}u \eta) \partial_{x}^{2}u\eta\mbox{d}x
-\int [H; \eta]\partial_{x}^{3}u\partial_{x}^{2}u\eta\mbox{d}x\nonumber\\
&=&\int  H\partial_{x}(\partial_{x}^{2}u \eta) \partial_{x}^{2}u\eta\mbox{d}x
-\int  H(\partial_{x}^{2}u \eta^{'}) \partial_{x}^{2}u\eta\mbox{d}x
-\int [H; \eta]\partial_{x}^{3}u\partial_{x}^{2}u\eta\mbox{d}x\nonumber\\
&=&A_{421}+A_{422}+A_{423}.
\end{eqnarray*}
Similar to the treatment of (\ref{1A421}), we write $A_{421}$ as
\begin{eqnarray*}
A_{421}
=\int  H\partial_{x}(\partial_{x}^{2}u \eta) \partial_{x}^{2}u\eta\mbox{d}x
=\int[D_{x}^{\frac{1}{2}}(\partial_{x}^{2}u \eta)]^{2}\mbox{d}x.
\end{eqnarray*}

The  Young inequality leads to
\begin{eqnarray}\label{2A422}
\int_{0}^{T}|A_{422}|\mbox{d}t
&\leq&\int_{0}^{T}\left|\int  H(\partial_{x}^{2}u \eta^{'}) \partial_{x}^{2}u\eta\right|\mbox{d}x\mbox{d}t\nonumber\\
&\leq& c_{0}\int_{0}^{T}\int(\partial_{x}^{2}u)^{2} (\eta^{'})^{2}\mbox{d}x\mbox{d}t
+c_{0}\int_{0}^{T}\int(\partial_{x}^{2}u \eta)^{2}\mbox{d}x\mbox{d}t.
\end{eqnarray}

Invoking (\ref{prp3}) and using (\ref{conclusion1}) with $(\epsilon/3, b+\epsilon)$ instead of
$(\epsilon, b)$ yield
\begin{eqnarray*}
&&\int_{0}^{T}\int(\partial_{x}^{2}u \eta)^{2}\mbox{d}x\mbox{d}t
+\int_{0}^{T}\int(\partial_{x}^{2}u)^{2} (\eta^{'})^{2}\mbox{d}x\mbox{d}t\nonumber\\
&&\leq \int_{0}^{T}\int(\partial_{x}^{2}u )^{2}\chi^{'}(x+vt; \epsilon, b)\mbox{d}x\mbox{d}t
+\int_{0}^{T}\int(\partial_{x}^{2}u)^{2}\chi^{'}(x+vt; \epsilon/3, b+\epsilon)\mbox{d}x\mbox{d}t\nonumber\\
&&\leq c_{0}.
\end{eqnarray*}
For the term $A_{423}$, (\ref{commutator}) leads to
\begin{eqnarray*}
A_{423}
&=&-\int [H; \eta]\partial_{x}^{3}u\partial_{x}^{2}u\eta\mbox{d}x\nonumber\\
&=&\int \partial_{x}[H; \eta]\partial_{x}^{3}u\partial_{x}u\eta\mbox{d}x
+\int [H; \eta]\partial_{x}^{3}u\partial_{x}u\eta^{'}\mbox{d}x\nonumber\\
&\leq& \|\partial_{x}[H; \eta]\partial_{x}^{3}u\|_{L^{2}}\|\partial_{x}u\eta\|_{L^{2}}
+\|[H; \eta]\partial_{x}^{3}u\|_{L^{2}}\|\partial_{x}u\eta^{'}\|_{L^{2}}\nonumber\\
&\leq& c\|u\|_{L^{2}}\|\partial_{x}u\eta\|_{L^{2}}
+c\|u\|_{L^{2}}\|\partial_{x}u\eta^{'}\|_{L^{2}}\nonumber\\
&\leq& c\|u_{0}\|_{L^{2}}^{2}+c\|\partial_{x}u\eta\|_{L^{2}}^{2}+c\|\partial_{x}u\eta^{'}\|_{L^{2}}^{2},
\end{eqnarray*}
which, after integration in time, can be controlled by using
similar method as that in (\ref{2A422}).

Accordingly, gathering the above information in (\ref{second}) and
invoking the Gronwall inequality, one derives

\begin{eqnarray}\label{conclusion2}
&&\sup_{t\in[0, T]}\int (\partial_{x}^{2}u)^{2}(x, t)\chi(x+vt; \epsilon, b)\mbox{d}x
+\int_{0}^{T}\int (\partial_{x}^{3}u)^{2}(x, t)\chi^{'}(x+vt; \epsilon, b)\mbox{d}x\mbox{d}t\nonumber\\
&&+\int_{0}^{T}\int [D_{x}^{\frac{1}{2}}(\partial_{x}^{2}u(x, t) \eta(x+vt; \epsilon, b))]^{2}\mbox{d}x\mbox{d}t
\leq c_{0}.
\end{eqnarray}

We prove the general case $ l \geq 2 $ by induction.
In details, we assume:
If $u_{0}$ satisfies (\ref{113}) then (\ref{114}) holds, that is to say
\begin{eqnarray}\label{l}
&&\sup_{t\in[0, T]}\int (\partial_{x}^{j}u)^{2}(x, t)\chi(x+vt; \epsilon, b)\mbox{d}x
+\int_{0}^{T}\int (\partial_{x}^{j+1}u)^{2}(x, t)\chi^{'}(x+vt; \epsilon, b)\mbox{d}x\mbox{d}t\nonumber\\
&&+\int_{0}^{T}\int [D_{x}^{\frac{1}{2}}(\partial_{x}^{j}u(x, t) \eta(x+vt; \epsilon, b))]^{2}\mbox{d}x\mbox{d}t
\leq c_{0}
\end{eqnarray}
for $j=1, 2, ..., l, $ $l\geq 2,$ and for any $\epsilon, b, v>0.$

Now we have that
\begin{eqnarray*}
u_{0}|_{(0, \infty)}\in H^{l+1}((0, \infty)).
\end{eqnarray*}
Thus from the previous step $(\ref{l})$ holds.
And formally, we have for $\epsilon, b, v>0$ the following identity
\begin{eqnarray}\label{third}
&&\frac{1}{2}\frac{\mbox{d}}{\mbox{d}t}\int(\partial_{x}^{l+1}u)^{2}(x, t)\chi(x+vt)\mbox{d}x
\underbrace{-\frac{1}{2}v\int(\partial_{x}^{l+1}u)^{2}(x, t)\chi^{'}(x+vt)\mbox{d}x}
_{\text{$A_{1}$}}\nonumber\\
&&+\frac{3}{2}\int(\partial_{x}^{l+2}u)^{2}(x, t)\chi^{'}(x+vt)\mbox{d}x
\underbrace{-\frac{1}{2}\int(\partial_{x}^{l+1}u)^{2}(x, t)\chi^{'''}(x+vt)\mbox{d}x}
_{\text{$A_{2}$}}\nonumber\\
&&\underbrace{+\int\partial_{x}^{l+1}(u\partial_{x}u)\partial_{x}^{l+1}u(x, t)\chi(x+vt)\mbox{d}x}
_{\text{$A_{3}$}}
\underbrace{-\int H\partial_{x}^{l+3}u\partial_{x}^{l+1}u(x, t)\chi(x+vt)\mbox{d}x}
_{\text{$A_{4}$}}=0.
\end{eqnarray}


Invoking (\ref{l}) with $j=l$, it holds that

\begin{eqnarray*}
\int_{0}^{T}|A_{1}(t)|\mbox{d}t
\leq |v|\int_{0}^{T}\int(\partial_{x}^{l+1}u)^{2}(x, t)\chi^{'}(x+vt)\mbox{d}x
\leq c_{0}.
\end{eqnarray*}
Using similar  method of treating (\ref{2A2}), we find
\begin{eqnarray*}
\int_{0}^{T}|A_{2}(t)|\mbox{d}t
&\leq&\int_{0}^{T}\int(\partial_{x}^{l+1}u)^{2}(x, t)|\chi^{'''}(x+vt; \epsilon, b)|
\mbox{d}x\mbox{d}t\nonumber\\
&\leq& \int_{0}^{T}\int(\partial_{x}^{l+1}u)^{2}(x, t)|\chi^{'}(x+vt; \epsilon/3, b+\epsilon)|\mbox{d}x\mbox{d}t\nonumber\\
&\leq& c_{0}.
\end{eqnarray*}

We estimate $A_{3}$ by considering two cases:
The first case is  when $l+1=3$ and the second is $l+1\geq 4.$

When $l+1=3$, we have after integration by parts
\begin{eqnarray*}
A_{3}(t)
&=&4\int\partial_{x}u(\partial_{x}^{3}u)^{2}\chi(x+vt)\mbox{d}x
+\int u \partial_{x}^{4}u\partial_{x}^{3}u\chi(x+vt)\mbox{d}x\nonumber\\
&+&3\int(\partial_{x}^{2}u)^{2}\partial_{x}^{3}u\chi(x+vt)\mbox{d}x\nonumber\\
&=&\frac{7}{2}\int\partial_{x}u(\partial_{x}^{3}u)^{2}\chi(x+vt)\mbox{d}x
-\frac{1}{2}\int u (\partial_{x}^{3}u)^{2}\chi^{'}(x+vt)\mbox{d}x\nonumber\\
&+&3\int(\partial_{x}^{2}u)^{2}\partial_{x}^{3}u\chi(x+vt)\mbox{d}x\nonumber\\
&=&A_{31}+A_{32}+A_{33}.
\end{eqnarray*}

Simple computation leads to
\begin{eqnarray*}
|A_{31}(t)|
<\|\partial_{x}u\|_{L^{\infty}}
\int(\partial_{x}^{3}u)^{2}\chi(x+vt)\mbox{d}x
\end{eqnarray*}
with the integral be the quantity to be estimated.

Employing (\ref{l}) with $j=l=2$, one deduces

\begin{eqnarray*}
\int_{0}^{T}|A_{32}(t)|\mbox{d}t
&\leq&\sup_{t\in[0, T]}\|u\|_{L^{\infty}}
\int_{0}^{T}\int(\partial_{x}^{3}u)^{2}\chi^{'}(x+vt)\mbox{d}x\mbox{d}t\nonumber\\
&\leq&\sup_{t\in[0, T]}\|u\|_{H^{3/4^{+}}}
\int_{0}^{T}\int(\partial_{x}^{3}u)^{2}\chi^{'}(x+vt)\mbox{d}x\mbox{d}t\nonumber\\
&\leq& c_{0}.
\end{eqnarray*}

Integration by parts leads to
\begin{eqnarray*}
A_{33}=
3\int(\partial_{x}^{2}u)^{2}\partial_{x}^{3}u\chi(x+vt)\mbox{d}x
=-\int(\partial_{x}^{2}u)^{3}\chi^{'}(x+vt)\mbox{d}x
\end{eqnarray*}
Using (\ref{id}), we have
\begin{eqnarray}\label{1A331}
|A_{33}|\leq
\|\partial_{x}^{2}u\chi^{'}(\cdot+vt; \epsilon, b )\|_{L^{\infty}}
\int(\partial_{x}^{2}u)^{2}\chi(x+vt; \epsilon/10, \epsilon/2 )\mbox{d}x
\end{eqnarray}
with the integral be bounded in $t\in (0, T]$ by a constant $c_{0}(\epsilon, b, v)$ resulting from
(\ref{l})(j=2).
Therefore, from the  boundedness of $\chi^{'}$
and the Sobolev inequality $\|f\|_{L^{\infty}}\leq \|f\|_{H^{1, 1}}$, one has
\begin{eqnarray}\label{1A332}
|A_{33}|
&\leq&
c\|\partial_{x}^{2}u\chi^{'}(\cdot+vt; \epsilon, b )\|_{L^{\infty}}^{2}+c\nonumber\\
&\leq& c\|(\partial_{x}^{2}u)^{2}\chi^{'}(\cdot+vt; \epsilon, b )\|_{L^{\infty}}+c\nonumber\\
&\leq& c\int|\partial_{x}[(\partial_{x}^{2}u)^{2}\chi^{'}(x+vt; \epsilon, b )]|\mbox{d}x
+c\nonumber\\
&\leq& c\int|\partial_{x}^{2}u\partial_{x}^{3}u\chi^{'}(x+vt; \epsilon, b )|\mbox{d}x
+c\int|\partial_{x}^{2}u\partial_{x}^{2}u\chi^{''}(x+vt; \epsilon, b )|\mbox{d}x
+c\nonumber\\
&\leq& c\int(\partial_{x}^{2}u)^{2}\chi^{'}(x+vt; \epsilon, b )|\mbox{d}x
+c\int(\partial_{x}^{3}u)^{2}\chi^{'}(x+vt; \epsilon, b )|\mbox{d}x\nonumber\\
&+&c\int|\partial_{x}^{2}u\partial_{x}^{2}u\chi^{'}(x+vt; \epsilon/3, b+\epsilon )|\mbox{d}x
+c,
\end{eqnarray}
where we have used (\ref{prp4}) with $j=2$.

Employing (\ref{l}) with $j=1, 2$ and integration in time, we obtain
\begin{eqnarray}\label{1A333}
\int_{0}^{T}|A_{33}|\mbox{d}t
\leq
c_{0}.
\end{eqnarray}

We turn our attention to the second case $l+1\geq 4$ in $A_{3}$.

By integration by parts, one derives
\begin{eqnarray*}
A_{3}
&=&d_{0}\int u(\partial_{x}^{l+1}u)^{2}\chi^{'}(x+vt)\mbox{d}x
+d_{1}\int \partial_{x}u(\partial_{x}^{l+1}u)^{2}\chi(x+vt)\mbox{d}x\nonumber\\
&+&d_{2}\int\partial_{x}^{2}u\partial_{x}^{l}u\partial_{x}^{l+1}u\chi(x+vt)\mbox{d}x
+\sum_{j=3}^{l-1}\int\partial_{x}^{j}u\partial_{x}^{l+2-j}u\partial_{x}^{l+1}u\chi(x+vt)\mbox{d}x\nonumber\\
&=&A_{3,0}+A_{3,1}+A_{3,2}+\sum_{j=3}^{l-1}A_{3, j}.
\end{eqnarray*}
Using (\ref{l}) with $j=l$ and the Sobolev embedding, one obtains
\begin{eqnarray*}
\int_{0}^{T}|A_{3, 0}|\mbox{d}t
&\leq&
\int_{0}^{T}\|u\|_{L^{\infty}}\int(\partial_{x}^{l+1}u)^{2}\chi^{'}(x+vt)\mbox{d}x\mbox{d}t\nonumber\\
&\leq& \sup_{0\leq t \leq T }\|u\|_{H^{3/4^{+}}}\int_{0}^{T}\int(\partial_{x}^{l+1}u)^{2}\chi^{'}(x+vt)\mbox{d}x\mbox{d}t
\leq c_{0}.
\end{eqnarray*}
Direct computation leads to
\begin{eqnarray*}
|A_{3, 1}|
\leq
\|\partial_{x}u\|_{L^{\infty}}\int(\partial_{x}^{l+1}u)^{2}\chi(x+vt)\mbox{d}x,
\end{eqnarray*}
which can be handled by the Gronwall inequality and (\ref{class})(ii).

To estimate $A_{3, 2}$ we follow the argument in the previous case.

Accordingly, we need to estimate $\sum_{j=3}^{l-1}A_{3, j}$ which only
appears when $l-1\geq 3.$

The Young inequality leads to
\begin{eqnarray*}
|A_{3, j}|
&\leq&
\frac{1}{2}\int(\partial_{x}^{j}u\partial_{x}^{l+2-j}u)^{2}\chi(x+vt)\mbox{d}x
+\frac{1}{2}\int(\partial_{x}^{l+1}u)^{2}\chi(x+vt)\mbox{d}x\nonumber\\
&=&A_{3, j, 1}+\frac{1}{2}\int(\partial_{x}^{l+1}u)^{2}\chi(x+vt)\mbox{d}x
\end{eqnarray*}
with the last integral be the quantity to be estimated.

To handle $A_{3, j, 1}$, one observes that $j, l+2-j\leq l-1$ and accordingly
\begin{eqnarray*}
|A_{3, j, 1}|
\leq
\|(\partial_{x}^{j}u)^{2}\chi(\cdot+vt; \epsilon/10,  \epsilon/2)\|_{L^{\infty}}
\int(\partial_{x}^{l+2-j}u)^{2}\chi(x+vt; \epsilon, b)\mbox{d}x
\end{eqnarray*}
with the last integral be bounded by (\ref{l}).
Moreover, Sobolev embedding yields
\begin{eqnarray*}
&&\|(\partial_{x}^{j}u)^{2}\chi(\cdot+vt; \epsilon/10,  \epsilon/2)\|_{L^{\infty}}\nonumber\\
&&\leq \|\partial_{x}[(\partial_{x}^{j}u)^{2}\chi(\cdot+vt; \epsilon/10,  \epsilon/2)]\|_{L^{1}}\nonumber\\
&&\leq \|\partial_{x}^{j}u\partial_{x}^{j+1}u\chi(\cdot+vt; \epsilon/10,  \epsilon/2)]\|_{L^{1}}
+\|\partial_{x}^{j}u\partial_{x}^{j}u\chi^{'}(\cdot+vt; \epsilon/10,  \epsilon/2)]\|_{L^{1}}\nonumber\\
&&\leq c\int(\partial_{x}^{j}u)^{2}\chi(x+vt; \epsilon/10,  \epsilon/2)\mbox{d}x
+c\int(\partial_{x}^{j+1}u)^{2}\chi(x+vt; \epsilon/10,  \epsilon/2)\mbox{d}x\nonumber\\
&&+c\int(\partial_{x}^{j}u)^{2}\chi^{'}(x+vt; \epsilon/10,  \epsilon/2)\mbox{d}x,
\end{eqnarray*}
which can be treated after integration in time by invoking (\ref{l}).

Finally, we estimate $A_{4}$.

After integration by parts, we find
\begin{eqnarray*}
A_{4}
&=&-\int H\partial_{x}^{l+3}u\partial_{x}^{l+1}u\chi(x+vt)\mbox{d}x\nonumber\\
&=&\int H\partial_{x}^{l+2}u\partial_{x}^{l+2}u\chi(x+vt)\mbox{d}x
+\int H\partial_{x}^{l+2}u\partial_{x}^{l+1}u\chi^{'}(x+vt)\mbox{d}x\nonumber\\
&=&A_{41}+A_{42}.
\end{eqnarray*}
Similar to (\ref{2A41}), we write $A_{41}$ as
\begin{eqnarray*}
A_{41}
&=&\int H\partial_{x}^{l+2}u\partial_{x}^{l+2}u\chi(x+vt)\mbox{d}x\nonumber\\
&=&-\int \partial_{x}^{l+2}u H(\partial_{x}^{l+2}u\chi(x+vt))\mbox{d}x\nonumber\\
\end{eqnarray*}
\begin{eqnarray*}
\quad\quad\quad\quad\quad\quad\quad\quad\quad\quad\quad&=&-\int \partial_{x}^{l+2}u H\partial_{x}^{l+2}u\chi(x+vt)\mbox{d}x
-\int \partial_{x}^{l+2}u [H; \chi]\partial_{x}^{l+2}u\mbox{d}x\nonumber\\
&=&-A_{41}-\int \partial_{x}^{l+2}u [H; \chi]\partial_{x}^{l+2}u\mbox{d}x.
\end{eqnarray*}
Consequently, there holds
\begin{eqnarray*}
A_{41}
&=&-\frac{1}{2}\int \partial_{x}^{l+2}u [H; \chi]\partial_{x}^{l+2}u\mbox{d}x\nonumber\\
&=&-\frac{1}{2}(-1)^{l+2}\int u \partial_{x}^{l+2}[H; \chi]\partial_{x}^{l+2}u\mbox{d}x\nonumber\\
&\leq& c\|u\|_{L^{2}}\|\partial_{x}^{l+2}[H; \chi]\partial_{x}^{l+2}u\|_{L^{2}}\nonumber\\
&\leq& c\|u\|_{L^{2}}^{2}=c\|u_{0}\|_{L^{2}}^{2}.
\end{eqnarray*}

Recall $\eta=\sqrt{\chi^{'}},$ therefore
\begin{eqnarray*}
A_{42}
&=&\int H\partial_{x}^{l+2}u\partial_{x}^{l+1}u\chi^{'}(x+vt)\mbox{d}x\nonumber\\
&=&\int H\partial_{x}^{l+2}u \eta \partial_{x}^{l+1}u\eta\mbox{d}x\nonumber\\
&=&\int  H(\partial_{x}^{l+2}u \eta) \partial_{x}^{l+1}u\eta\mbox{d}x
-\int [H; \eta]\partial_{x}^{l+2}u\partial_{x}^{l+1}u\eta\mbox{d}x\nonumber\\
&=&\int  H\partial_{x}(\partial_{x}^{l+1}u \eta) \partial_{x}^{l+1}u\eta\mbox{d}x
-\int  H(\partial_{x}^{l+1}u \eta^{'}) \partial_{x}^{l+1}u\eta\mbox{d}x\nonumber\\
&-&\int [H; \eta]\partial_{x}^{l+2}u\partial_{x}^{l+1}u\eta\mbox{d}x\nonumber\\
&=&A_{421}+A_{422}+A_{423}.
\end{eqnarray*}
For the term $A_{421}$, one has
\begin{eqnarray*}
A_{421}=\int  H\partial_{x}(\partial_{x}^{l+1}u \eta) \partial_{x}^{l+1}u\eta\mbox{d}x
=\int[D_{x}^{\frac{1}{2}}(\partial_{x}^{l+1}u \eta)]^{2}\mbox{d}x.
\end{eqnarray*}

The H\"{o}lder and Young inequality yield
\begin{eqnarray}\label{3A422}
\int_{0}^{T}|A_{422}|\mbox{d}t
&\leq&\int_{0}^{T}\left|\int  H(\partial_{x}^{l+1}u \eta^{'}) \partial_{x}^{l+1}u\eta\right|\mbox{d}x\mbox{d}t\nonumber\\
&\leq& c\int_{0}^{T}\int(\partial_{x}^{l+1}u)^{2} (\eta^{'})^{2}\mbox{d}x\mbox{d}t
+c\int_{0}^{T}\int(\partial_{x}^{l+1}u \eta)^{2}\mbox{d}x\mbox{d}t.
\end{eqnarray}
Thus, we can handle this term by using a similar method
as that in (\ref{2A422}).

Invoking (\ref{commutator}), one finds
\begin{eqnarray*}
A_{423}
&=&-\int [H; \eta]\partial_{x}^{l+2}u\partial_{x}^{l+1}u\eta\mbox{d}x\nonumber\\
&=&\int \partial_{x}[H; \eta]\partial_{x}^{l+2}u\partial_{x}^{l}u\eta\mbox{d}x
+\int [H; \eta]\partial_{x}^{l+2}u\partial_{x}^{l}u\eta^{'}\mbox{d}x\nonumber\\
&\leq& \|\partial_{x}[H; \eta]\partial_{x}^{l+2}u\|_{L^{2}}\|\partial_{x}^{l}u\eta\|_{L^{2}}
+\|[H; \eta]\partial_{x}^{l+2}u\|_{L^{2}}\|\partial_{x}^{l}u\eta^{'}\|_{L^{2}}\nonumber\\
&\leq& c\|u\|_{L^{2}}\|\partial_{x}^{l}u\eta\|_{L^{2}}
+c\|u\|_{L^{2}}\|\partial_{x}^{l}u\eta^{'}\|_{L^{2}}\nonumber\\
&\leq& c\|u_{0}\|_{L^{2}}^{2}+c\|\partial_{x}^{l}u\eta\|_{L^{2}}^{2}+c\|\partial_{x}^{l}u\eta^{'}\|_{L^{2}}^{2},
\end{eqnarray*}
which can also be controlled by using a similar way as that in (\ref{2A422}).

As a consequence, substituting  the above information into (\ref{third}) and
employing the Gronwall inequality, one deduces

\begin{eqnarray}\label{conclusion3}
&&\sup_{t\in[0, T]}\int (\partial_{x}^{l+1}u)^{2}\chi(x+vt; \epsilon, b)\mbox{d}x
+\int_{0}^{T}\int (\partial_{x}^{l+2}u)^{2}\chi^{'}(x+vt; \epsilon, b)\mbox{d}x\mbox{d}t\nonumber\\
&&+\int_{0}^{T}\int [D_{x}^{\frac{1}{2}}(\partial_{x}^{l+1}u \eta)]^{2}\mbox{d}x\mbox{d}t
\leq c_{0}.
\end{eqnarray}
This close our induction.

To justify the previous formal computations we
refer the reader to \cite{ilp1} and we omit the details here.

\section{Proof of Theorem \ref{persistence}}
In this section, we prove Theorem \ref{persistence}.

We first prove (\ref{118}) for any $n\in\mathbb{Z}^{+}.$

Note that $x_{+}^{n}u_{0}\in L^{2}(\mathbb{R})$
implies $\chi_{n}(x; \epsilon, b)u_{0} \in L^{2}(\mathbb{R}).$

Multiplying equation (\ref{benjamin}) with $u(x, t)\chi_{n}(x+vt; \epsilon, b)$ and integrating,
we obtain
\begin{eqnarray}\label{zero}
&&\frac{1}{2}\frac{\mbox{d}}{\mbox{d}t}\int u^{2}(x, t)\chi_{n}(x+vt; \epsilon, b)\mbox{d}x
\underbrace{-\frac{1}{2}v\int u^{2}(x, t)\chi^{'}_{n}(x+vt; \epsilon, b)\mbox{d}x}
_{\text{$A_{1}$}}\nonumber\\
&&+\frac{3}{2}\int(\partial_{x}u)^{2}(x, t)\chi^{'}_{n}(x+vt; \epsilon, b)\mbox{d}x
\underbrace{-\frac{1}{2}\int  u^{2}(x, t)\chi^{'''}_{n}(x+vt; \epsilon, b)\mbox{d}x}
_{\text{$A_{2}$}}\nonumber\\
&&\underbrace{+\int u\partial_{x}uu(x, t)\chi_{n}(x+vt; \epsilon, b)\mbox{d}x}
_{\text{$A_{3}$}}
\underbrace{-\int H\partial_{x}^{2}uu(x, t)\chi_{n}(x+vt; \epsilon, b)\mbox{d}x}
_{\text{$A_{4}$}}=0.
\end{eqnarray}

Employing (\ref{prp8}) with $j=1$, one easily deduces
\begin{eqnarray*}\label{46}
|A_{1}(t)|
&\leq& |v|\int u^{2}\chi_{n}(x+vt; \epsilon, b)\mbox{d}x
+|v|c(n, b)\int u^{2}\mbox{d}x\nonumber\\
&\leq& |v|\int u^{2}\chi_{n}(x+vt; \epsilon, b)\mbox{d}x
+|v|c(n, b)\|u_{0}\|_{L^{2}}^{2}.
\end{eqnarray*}
Again, invoking  (\ref{prp8}) with $j=3,$ we derive
\begin{eqnarray*}\label{45}
|A_{2}(t)|
&\leq& c(n, b)\int u^{2}\mbox{d}x
+\int u^{2}\chi_{n}(x+vt; \epsilon, b)\mbox{d}x\nonumber\\
&\leq& c(n, b)\|u_{0}\|_{L^{2}}^{2}
+\int u^{2}\chi_{n}(x+vt; \epsilon, b)\mbox{d}x.
\end{eqnarray*}

 We estimate $A_{3}$ as
\begin{eqnarray*}\label{47}
|A_{3}(t)|
\leq \|\partial_{x}u\|_{L^{\infty}}\int u^{2}\chi_{n}(x+vt; \epsilon, b)\mbox{d}x.
\end{eqnarray*}
For $A_{4}$, integration by parts produces
\begin{eqnarray*}\label{48}
A_{4}
&=&-\int H\partial_{x}^{2}uu\chi_{n}(x+vt; \epsilon, b)\mbox{d}x\nonumber\\
&=&\int H\partial_{x}u\partial_{x}u\chi_{n}(x+vt; \epsilon, b)\mbox{d}x
+\int H\partial_{x}uu\chi_{n}^{'}(x+vt; \epsilon, b)\mbox{d}x\nonumber\\
&=&A_{41}+A_{42}.
\end{eqnarray*}
Reasoning as (\ref{2A41}), it holds that
\begin{eqnarray*}
A_{41}
&=&\int H\partial_{x}u\partial_{x}u\chi_{n}(x+vt; \epsilon, b)\mbox{d}x\nonumber\\
&=&-\int \partial_{x}u H(\partial_{x}u\chi_{n}(x+vt; \epsilon, b))\mbox{d}x\nonumber\\
&=&-\int \partial_{x}u H\partial_{x}u\chi_{n}(x+vt; \epsilon, b)\mbox{d}x
-\int \partial_{x}u [H; \chi_{n}]\partial_{x}u\mbox{d}x\nonumber\\
&=&-A_{41}-\int \partial_{x}u [H; \chi_{n}]\partial_{x}u\mbox{d}x.
\end{eqnarray*}
As a result, we find
\begin{eqnarray*}
A_{41}
&=&-\frac{1}{2}\int \partial_{x}u [H; \chi_{n}]\partial_{x}u\mbox{d}x\nonumber\\
&=&\frac{1}{2}\int u \partial_{x}[H; \chi_{n}]\partial_{x}u\mbox{d}x\nonumber\\
&\leq& c\|u\|_{L^{2}}\|\partial_{x}[H; \chi_{n}]\partial_{x}u\|_{L^{2}}\nonumber\\
&\leq& c\|u\|_{L^{2}}^{2}=c\|u_{0}\|_{L^{2}}^{2}.
\end{eqnarray*}
Integration by parts  and (\ref{eta}) lead to
\begin{eqnarray}\label{4A42}
A_{42}
&=&\int H\partial_{x}uu\chi_{n}^{'}(x+vt)\mbox{d}x\nonumber\\
&=&\int H\partial_{x}u \eta_{n}u\eta_{n}\mbox{d}x\nonumber\\
&=&\int  H(\partial_{x}u \eta_{n}) u\eta_{n}\mbox{d}x
-\int [H; \eta_{n}]\partial_{x}uu\eta_{n}\mbox{d}x\nonumber\\
&=&\int  H\partial_{x}(u \eta_{n}) u\eta_{n}\mbox{d}x
-\int  H(u \eta_{n}^{'}) u\eta_{n}\mbox{d}x
-\int [H; \eta_{n}]\partial_{x}uu\eta_{n}\mbox{d}x\nonumber\\
&=&A_{421}+A_{422}+A_{423}.
\end{eqnarray}
Again, using the  Plancherel theorem, we write $A_{421}$ as
\begin{eqnarray*}
A_{421}=\int  H\partial_{x}(u \eta_{n}) u\eta_{n}\mbox{d}x
=\int[D_{x}^{\frac{1}{2}}(u \eta_{n})]^{2}\mbox{d}x.
\end{eqnarray*}

For the term $A_{422}$, note that because of the factor $x^{n},$
the support of $\chi_{n}^{'}$ is not $[\epsilon, b]$ at all for the general case.
As a consequence,  $\eta_{n}$
and $\eta_{n}^{'}$ may be unbounded.
However, we notice that (\ref{prp10}) and (\ref{prp11})
provide us with a relation between $\chi_{n} $ and $\chi_{n-1} $,
therefore, we could use induction to close our proof.

Let us first consider the case $n=0$ and thus $\eta_{n}=\eta$ .

At this point, we derive
\begin{eqnarray*}\label{4A422}
A_{422}=
-\int  H(u \eta^{'}) u\eta\mbox{d}x
\leq c\|u \eta^{'}\|_{L^{2}}\|u \eta\|_{L^{2}}
\leq c\|u \|_{L^{2}}^{2}
\leq c\|u_{0} \|_{L^{2}}^{2}.
\end{eqnarray*}

Invoking (\ref{commutator}), we find
\begin{eqnarray*}
A_{423}
=-\int [H; \eta]\partial_{x}uu\eta\mbox{d}x
\leq \|[H; \eta]\partial_{x}u\|_{L^{2}}\|u\eta\|_{L^{2}}
\leq c\|u \|_{L^{2}}^{2}
\leq c\|u_{0} \|_{L^{2}}^{2}.
\end{eqnarray*}
Hence, we obtain the following inequality when $n=0:$
\begin{eqnarray*}\label{n0}
&&\sup_{t\in[0, T]}\int u^{2}(x, t)\chi(x+vt; \epsilon, b)\mbox{d}x
+\int_{0}^{T}\int (\partial_{x}u)^{2}(x, t)\chi^{'}(x+vt; \epsilon, b)\mbox{d}x\mbox{d}t\nonumber\\
&&+\int_{0}^{T}\int [D_{x}^{\frac{1}{2}}(u(x, t) \eta(x+vt; \epsilon, b))]^{2}\mbox{d}x\mbox{d}t
\leq c_{2}.
\end{eqnarray*}
Let us assume the case $n\geq 0$ holds, i.e.,
\begin{eqnarray}\label{nn}
&&\sup_{t\in[0, T]}\int u^{2}(x, t)\chi_{n}(x+vt; \epsilon, b)\mbox{d}x
+\int_{0}^{T}\int (\partial_{x}u)^{2}(x, t)\chi_{n}^{'}(x+vt; \epsilon, b)\mbox{d}x\mbox{d}t\nonumber\\
&&+\int_{0}^{T}\int [D_{x}^{\frac{1}{2}}(u(x, t) \eta_{n}(x+vt; \epsilon, b))]^{2}\mbox{d}x\mbox{d}t
\leq c_{2}.
\end{eqnarray}
We shall prove the case $n+1$.

We only need to treat the terms $A_{422}$ and $A_{423}$ with $n+1$
instead of $n$.

Employing (\ref{prp10}) and (\ref{prp11}), we find
\begin{eqnarray*}\label{4A4222}
A_{422}
&=&
-\int  H(u \eta_{n+1}^{'}) u\eta_{n+1}\mbox{d}x\nonumber\\
&\leq& c_{2}\|u \eta_{n+1}^{'}\|_{L^{2}}\|u \eta_{n+1}\|_{L^{2}}\nonumber\\
&\leq& c_{2}\int u^{2} (\eta_{n+1}^{'})^{2}\mbox{d}x
+c_{2}\int u^{2} (\eta_{n+1})^{2}\mbox{d}x\nonumber\\
&\leq& c_{2}\int u^{2} \chi_{n}(x+vt; \epsilon/3, b+\epsilon)\mbox{d}x,
\end{eqnarray*}
which can be handled by using (\ref{nn}) with $(\epsilon/3, b+\epsilon)$
instead of $(\epsilon, b).$

The term $A_{423}$ can be controlled similarly.

Thus we completes the proof of (\ref{118}).
And for convenience, we view (\ref{nn}) as a conclusion in the following of this paper.

Next, we prove (\ref{119}).

We first prove the case $n=1.$

From (\ref{nn}) with $n=1$ and (\ref{prp2}), it follows that for any $\delta> 0$
there exists  $\hat{t}\in (0, \delta)$ such that
\begin{eqnarray*}\label{410}
\int(\partial_{x}u)^{2}(x, \hat{t})\chi(x; \epsilon, b)\mbox{d}x
<\infty.
\end{eqnarray*}

A smooth solution $u$ to the  IVP (\ref{benjamin})
satisfies the following identity:
\begin{eqnarray}
&&\frac{1}{2}\frac{\mbox{d}}{\mbox{d}t}\int(\partial_{x}u)^{2}(x, t)\chi(x+vt; \epsilon, b)\mbox{d}x
\underbrace{-\frac{1}{2}v\int(\partial_{x}u)^{2}(x, t)\chi^{'}(x+vt; \epsilon, b)\mbox{d}x}
_{\text{$A_{1}$}}\nonumber
\end{eqnarray}
\begin{eqnarray}\label{411}
&&+\frac{3}{2}\int(\partial_{x}^{2}u)^{2}(x, t)\chi^{'}(x+vt; \epsilon, b)\mbox{d}x
\underbrace{-\frac{1}{2}\int(\partial_{x}u)^{2}(x, t)\chi^{'''}(x+vt; \epsilon, b)\mbox{d}x}
_{\text{$A_{2}$}}\nonumber\\
&&\underbrace{+\int\partial_{x}(u\partial_{x}u)\partial_{x}u(x, t)\chi(x+vt; \epsilon, b)\mbox{d}x}
_{\text{$A_{3}$}}
\underbrace{-\int H\partial_{x}^{3}u(x, t)\partial_{x}u\chi(x+vt; \epsilon, b)\mbox{d}x}
_{\text{$A_{4}$}}=0.
\end{eqnarray}

Using (\ref{nn}) with $n=0$, one obtains
\begin{eqnarray*}\label{412}
\int_{\hat{t}}^{T}|A_{1}(t)|\mbox{d}t
\leq c_{3}.
\end{eqnarray*}

Again, using (\ref{nn}) with $n=0$
and $(\epsilon/3, b+\epsilon)$ instead of $(\epsilon, b)$,
there holds
\begin{eqnarray*}\label{414}
\int_{\hat{t}}^{T}|A_{2}(t)|\mbox{d}t
\leq c_{3}.
\end{eqnarray*}

For the term $A_{3}$, integration by parts leads to

\begin{eqnarray*}\label{415}
A_{3}(t)
&=&\frac{1}{2}\int(\partial_{x}u)^{3}\chi(x+vt; \epsilon, b)\mbox{d}x
-\frac{1}{2}\int u \partial_{x}u\partial_{x}u\chi^{'}(x+vt; \epsilon, b)\mbox{d}x\nonumber\\
&\leq& \|\partial_{x}u\|_{L^{\infty}}\int(\partial_{x}u)^{2}\chi(x+vt; \epsilon, b)\mbox{d}x
+\|u\|_{L^{\infty}}\int(\partial_{x}u)^{2}\chi^{'}(x+vt; \epsilon, b)\mbox{d}x\nonumber\\
&=&A_{31}+A_{32}.
\end{eqnarray*}
By Sobolev embedding theorem, one obtains
\begin{eqnarray*}
\int_{\hat{t}}^{T}|A_{32}(t)|\mbox{d}t
<\sup_{t\in[\hat{t}, T]}\|u\|_{H^{3/4^{+}}}
\int_{\hat{t}}^{T}\int(\partial_{x}u)^{2}\chi^{'}(x+vt; \epsilon, b)\mbox{d}x\mbox{d}t.
\end{eqnarray*}
The term $A_{31}$ will be controlled by using
(\ref{class})(ii) and the Gronwall inequality.

The term $A_{4}$ can be estimated as in the proof of
Theorem \ref{regularity}, we omit it.
%

Substituting  the above information in $(\ref{411})$, using
Gronwall inequality and (\ref{class})(ii), one obtains

\begin{eqnarray}\label{conclusion5}
&&\sup_{t\in[\hat{t}, T]}\int (\partial_{x}u)^{2}(x, t)\chi(x+vt; \epsilon, b)\mbox{d}x
+\int_{\hat{t}}^{T}\int (\partial_{x}^{2}u)^{2}(x, t)\chi^{'}(x+vt; \epsilon, b)\mbox{d}x\mbox{d}t\nonumber\\
&&+\int_{\hat{t}}^{T}\int [D_{x}^{\frac{1}{2}}(\partial_{x}u(x, t) \eta(x+vt; \epsilon, b))]^{2}\mbox{d}x\mbox{d}t
\leq c_{3}.
\end{eqnarray}

Next, we turn to the case $n=2$ in the proof of (\ref{119}).

Since  $x_{+}u_{0}\in L^{2}(\mathbb{R})$, using (\ref{nn}) with $n=2$, one finds
\begin{eqnarray}\label{419}
&&\sup_{t\in[0, T]}\int u^{2}(x, t)\chi_{2}(x+vt; \epsilon, b)\mbox{d}x
+\int_{0}^{T}\int (\partial_{x}u)^{2}(x, t)\chi_{2}^{'}(x+vt; \epsilon, b)\mbox{d}x\mbox{d}t\nonumber\\
&&+\int_{0}^{T}\int [D_{x}^{\frac{1}{2}}(u(x, t) \eta_{2}(x+vt; \epsilon, b))]^{2}\mbox{d}x\mbox{d}t
\leq c_{2}.
\end{eqnarray}

Using  (\ref{419}) and (\ref{prp7}), we derive that for any $\delta> 0$
 there exists $\hat{t}\in (0, \delta)$ such that
 \begin{eqnarray*}\label{421}
\int (\partial_{x}u)^{2}(x, \hat{t})\chi_{1}(x; \epsilon, b)\mbox{d}x
<\infty.
\end{eqnarray*}

Consider the following identity
\begin{eqnarray}\label{z411}
&&\frac{1}{2}\frac{\mbox{d}}{\mbox{d}t}\int(\partial_{x}u)^{2}\chi_{1}(x+vt; \epsilon, b)\mbox{d}x
\underbrace{-\frac{1}{2}v\int(\partial_{x}u)^{2}(x, t)\chi_{1}^{'}(x+vt; \epsilon, b)\mbox{d}x}
_{\text{$A_{1}$}}\nonumber\\
&&+\frac{3}{2}\int(\partial_{x}^{2}u)^{2}(x, t)\chi_{1}^{'}(x+vt; \epsilon, b)\mbox{d}x
\underbrace{-\frac{1}{2}\int(\partial_{x}u)^{2}(x, t)\chi_{1}^{'''}(x+vt; \epsilon, b)\mbox{d}x}
_{\text{$A_{2}$}}\\
&&\underbrace{+\int\partial_{x}(u\partial_{x}u)(x, t)\partial_{x}u\chi_{1}(x+vt; \epsilon, b)\mbox{d}x}
_{\text{$A_{3}$}}
\underbrace{-\int H\partial_{x}^{3}u\partial_{x}u(x, t)\chi_{1}(x+vt; \epsilon, b)\mbox{d}x}
_{\text{$A_{4}$}}=0.\nonumber
\end{eqnarray}

Invoking  (\ref{nn}) with  $n=1$, it holds that
\begin{eqnarray*}\label{423}
\int_{\hat{t}}^{T}|A_{1}(t)|\mbox{d}t
\leq |v|\int_{\hat{t}}^{T}\int(\partial_{x}u)^{2}(x, t)\chi_{1}^{'}(x+vt; \epsilon, b)\mbox{d}x\mbox{d}t
\leq c_{3}.
\end{eqnarray*}

For the term $A_{2}$, employing the fact that $\chi_{1}^{'''}$
is supported in $[\epsilon, b]$, we deduce
\begin{eqnarray*}\label{426}
\int_{\hat{t}}^{T}|A_{2}(t)|\mbox{d}t
\leq c_{3},
\end{eqnarray*}
where we have used (\ref{class})(iii) with $r=0$.

Concerning  the term $A_{3}$, integration by parts leads to

\begin{eqnarray*}\label{w415}
A_{3}(t)
&=&\frac{1}{2}\int(\partial_{x}u)^{3}\chi_{1}(x+vt; \epsilon, b)\mbox{d}x
-\frac{1}{2}\int u \partial_{x}u\partial_{x}u\chi_{1}^{'}(x+vt; \epsilon, b)\mbox{d}x\nonumber\\
&\leq& \|\partial_{x}u\|_{L^{\infty}}\int(\partial_{x}u)^{2}\chi_{1}(x+vt; \epsilon, b)\mbox{d}x
+\|u\|_{L^{\infty}}\int(\partial_{x}u)^{2}\chi_{1}^{'}(x+vt; \epsilon, b)\mbox{d}x\nonumber\\
&=&A_{31}+A_{32},
\end{eqnarray*}
which can treated as in the proof of Theorem \ref{regularity}.

Finally, we control $A_{4}$.

Integration by parts yields
\begin{eqnarray*}
A_{4}
&=&-\int H\partial_{x}^{3}u\partial_{x}u\chi_{1}(x+vt; \epsilon, b)\mbox{d}x\nonumber\\
&=&\int H\partial_{x}^{2}u\partial_{x}^{2}u\chi_{1}(x+vt; \epsilon, b)\mbox{d}x
+\int H\partial_{x}^{2}u\partial_{x}u\chi_{1}^{'}(x+vt; \epsilon, b)\mbox{d}x\nonumber\\
&=&A_{41}+A_{42}.
\end{eqnarray*}
The term $A_{41}$ can be handled by using integration
by parts, the skew symmetry of the Hilbert transform and
the commutator estimate (\ref{commutator}), we omit it.

Now, we focus on the term $A_{42}.$
Recall $(\eta_{1})^{2}=\chi_{1}^{'},$ one has
\begin{eqnarray*}
A_{42}
&=&\int H\partial_{x}^{2}u\partial_{x}u\chi_{1}^{'}(x+vt)\mbox{d}x\nonumber\\
&=&\int H\partial_{x}^{2}u \eta_{1} \partial_{x}u\eta_{1}\mbox{d}x\nonumber
\end{eqnarray*}
\begin{eqnarray*}
&=&\int  H(\partial_{x}^{2}u \eta_{1}) \partial_{x}u\eta_{1}\mbox{d}x
-\int [H; \eta_{1}]\partial_{x}^{2}u\partial_{x}u\eta_{1}\mbox{d}x\nonumber\\
&=&\int  H\partial_{x}(\partial_{x}u \eta_{1}) \partial_{x}u\eta_{1}\mbox{d}x
-\int  H(\partial_{x}u \eta_{1}^{'}) \partial_{x}u\eta_{1}\mbox{d}x
-\int [H; \eta_{1}]\partial_{x}^{2}u\partial_{x}u\eta_{1}\mbox{d}x\nonumber\\
&=&A_{421}+A_{422}+A_{423}.
\end{eqnarray*}
For the term $A_{421}$, we have
\begin{eqnarray*}
A_{421}=\int  H\partial_{x}(\partial_{x}u \eta_{1}) \partial_{x}u\eta_{1}\mbox{d}x
=\int[D_{x}^{\frac{1}{2}}(\partial_{x}u \eta_{1})]^{2}\mbox{d}x.
\end{eqnarray*}

The Young inequality leads to
\begin{eqnarray*}\label{5A422}
\int_{\hat{t}}^{T}|A_{422}|\mbox{d}t
&=&\int_{\hat{t}}^{T}\left|\int  H(\partial_{x}u \eta_{1}^{'}) \partial_{x}u\eta_{1}\mbox{d}x\right|\mbox{d}t\nonumber\\
&\leq& c_{3}\int_{\hat{t}}^{T}\int(\partial_{x}u)^{2} (\eta_{1}^{'})^{2}\mbox{d}x\mbox{d}t
+c_{3}\int_{\hat{t}}^{T}\int(\partial_{x}u)^{2} \eta_{1}^{2}\mbox{d}x\mbox{d}t.
\end{eqnarray*}
Note that $\eta_{1}$ in unbounded in support of $\chi_{1}^{'}$.
However, invoking (\ref{nn}) with $n=1,$ we find
\begin{eqnarray*}
\int_{\hat{t}}^{T}\int(\partial_{x}u)^{2} \eta_{1}^{2}\mbox{d}x\mbox{d}t
\leq c_{3}.
\end{eqnarray*}
Now, simple computation yields $\chi(x; \epsilon, b)\leq \chi_{1}^{'}(x; \epsilon, b)$.
This fact combining with (\ref{prp11}) and (\ref{nn}) with $(\epsilon/3, b+\epsilon)$ instead of
$(\epsilon, b)$ permits us to conclude
\begin{eqnarray*}
\int_{\hat{t}}^{T}\int(\partial_{x}u)^{2} (\eta_{1}^{'})^{2}\mbox{d}x\mbox{d}t
\leq c_{3}.
\end{eqnarray*}
Thus we have controlled $A_{422}$ after integration in time.
The term $A_{423}$ can be handled by using the above method and (\ref{commutator}),
we omit it.

As a result, we conclude after invoking the Gronwall inequality that
\begin{eqnarray}\label{conclusion6}
&&\sup_{t\in[\hat{t}, T]}\int (\partial_{x}u)^{2}(x, t)\chi_{1}(x+vt; \epsilon, b)\mbox{d}x
+\int_{\hat{t}}^{T}\int (\partial_{x}^{2}u)^{2}(x, t)\chi_{1}^{'}(x+vt; \epsilon, b)\mbox{d}x\mbox{d}t\nonumber\\
&&+\int_{\hat{t}}^{T}\int [D_{x}^{\frac{1}{2}}(\partial_{x}u(x, t) \eta_{1}(x+vt; \epsilon, b))]^{2}\mbox{d}x\mbox{d}t
\leq c_{3}.
\end{eqnarray}

By (\ref{conclusion6}) for any $\delta> 0$ there exists
$\hat{\hat{t}}\in (\hat{t}, \delta)$ such that
\begin{eqnarray*}
\int (\partial_{x}^{2}u)^{2}(x, \hat{\hat{t}})\chi_{1}^{'}(x; \epsilon, b)\mbox{d}x< \infty,
\end{eqnarray*}
this produces
\begin{eqnarray*}
\int (\partial_{x}^{2}u)^{2}(x, \hat{\hat{t}})\chi(x; \epsilon, b)\mbox{d}x< \infty.
\end{eqnarray*}

Hence, the  result of  propagation of regularity (\ref{conclusion2}) yields :
\begin{eqnarray*}\label{431}
&&\sup_{t\in[\delta, T]}\int (\partial_{x}^{2}u)^{2}(x, t)\chi(x+vt; \epsilon, b)\mbox{d}x
+\int_{\delta}^{T}\int (\partial_{x}^{3}u)^{2}(x, t)\chi^{'}(x+vt; \epsilon, b)\mbox{d}x\mbox{d}t\nonumber\\
&&+\int_{\delta}^{T}\int [D_{x}^{\frac{1}{2}}(\partial_{x}^{2}u(x, t) \eta(x+vt; \epsilon, b))]^{2}\mbox{d}x\mbox{d}t
\leq c_{3}.
\end{eqnarray*}
This completes the proof of the case $n=2$.

For the general case, we use induction.

Given $(m, l)\in \mathbb{Z}^{+}\times \mathbb{Z}^{+}$
we say that

\begin{eqnarray}\label{432}
(m, l)>(\hat{m}, \hat{l})\Leftrightarrow
\begin{cases}
(i)\quad m>\hat{m}\\
or\\
(ii) \quad m=\hat{m}\quad and\quad l>\hat{l}.
\end{cases}
\end{eqnarray}
Similarly, we say that $(m, l)\geq(\hat{m}, \hat{l})$ if $(ii)$
in the right hand side of (\ref{432}) holds with $\geq$ instead of $>$.

The general case $(m, l)$ reads:

For any $\epsilon, b, v>0$
\begin{eqnarray}\label{433}
&&\sup_{t \in [\delta, T]}\int (\partial_{x}^{l}u)^{2}(x, t)\chi_{m}(x+vt; \epsilon, b)\mbox{d}x
+\int_{\delta}^{T}\int (\partial_{x}^{l+1}u)^{2}(x, t)\chi_{m}^{'}(x+vt; \epsilon, b)\mbox{d}x\mbox{d}t\nonumber\\
&&+\int_{\delta}^{T}\int [D_{x}^{\frac{1}{2}}(\partial_{x}^{l}u(x, t) \eta_{m}(x+vt; \epsilon, b))]^{2}\mbox{d}x\mbox{d}t
\leq c_{3}.
\end{eqnarray}

Notice that we have already proved the following cases:
\begin{enumerate}
  \item (0,1) and (1, 0)
  \item (0, 2) (1, 1) and  (2, 0)
  \item Under the hypothesis $x_{+}^{n/2}u_{0}\in L^{2}(\mathbb{R})$, we proved (\ref{nn}),
   i.e. (n, 0) for all $n\in \mathbb{Z}^{+}$
  \item By Theorem \ref{regularity} (propagation of regularity):
  If (\ref{433}) holds with $(m, l)=(1, l)$ $(\delta/2 $ instead of $\delta)$,
  then there exists $\hat{t}\in (\delta/2 , \delta)$ such that
  \begin{eqnarray*}
  \int (\partial_{x}^{l+1}u)^{2}\chi_{1}^{'}(x+v\hat{t}; \epsilon, b)\mbox{d}x
  < \infty
  \end{eqnarray*}
  which implies that
  \begin{eqnarray*}
  \int (\partial_{x}^{l+1}u)^{2}\chi(x+v\hat{t}; \epsilon, b)\mbox{d}x
  < \infty.
  \end{eqnarray*}
\end{enumerate}

By the propagation of regularity (Theorem \ref{regularity}), one has the result
(\ref{433}) with $(m, l)=(0, l+1)$, that is, $(1, l)$ implies $(0, l+1)$ for any
$l\in \mathbb{Z}^{+}$.

Now we assume
(\ref{433}) holds for $(m, k)$ such that

\begin{eqnarray*}\label{434}
 \begin{cases}
(a) \quad (m, k) \leq (n-j, j)\quad for \quad some \quad j=0, 1, 2,..., n\\
and\\
(b)\quad (m, k) = (n+1, 0), (n, 1), ..., (n+1-l, l)\quad for \quad some \quad l\leq n.
  \end{cases}
  \end{eqnarray*}

We need to prove the case $(n+1-(l+1), l+1)=(n-l, l+1)$.
From (4) above, since $(1, l)$ implies $(0, l+1)$, this case is already true for $l=n.$
Thus it suffices to consider $l\leq n-1.$

From the previous step $(n-l+1, l)$ we have that for any $\delta^{'}, v, \epsilon> 0,$

\begin{eqnarray}\label{435}
&&\sup_{t\in[\delta^{'}, T]}\int (\partial_{x}^{l}u)^{2}(x, t)\chi_{n+1-l}(x+vt; \epsilon, b)\mbox{d}x
+\int_{\delta^{'}}^{T}\int (\partial_{x}^{l+1}u)^{2}(x, t)\chi_{n+1-l}^{'}(x+vt; \epsilon, b)\mbox{d}x\mbox{d}t\nonumber\\
&&+\int_{\delta^{'}}^{T}\int [D_{x}^{\frac{1}{2}}(\partial_{x}^{l}u(x, t) \eta_{n+1-l}(x+vt; \epsilon, b))]^{2}\mbox{d}x\mbox{d}t
\leq c_{3}.
\end{eqnarray}

Simple computation yields
\begin{eqnarray*}\label{436}
\chi_{n+1-l}^{'}(x; \epsilon, b)\geq c\chi_{n-l}(x; \epsilon, b).
\end{eqnarray*}
According to (\ref{435}), there exists $\hat{t}\in (\delta^{'}, 2\delta^{'})$
such that
\begin{eqnarray*}\label{437}
\int (\partial_{x}^{l+1}u)^{2}(x, \hat{t})\chi_{n-l}(x+v\hat{t}; \epsilon, b)\mbox{d}x
< \infty.
\end{eqnarray*}
For smooth solution of equation (\ref{benjamin}), consider
\begin{eqnarray}\label{438}
&&\frac{1}{2}\frac{\mbox{d}}{\mbox{d}t}\int(\partial_{x}^{l+1}u)^{2}\chi_{n-l}(x+vt; \epsilon, b)\mbox{d}x
\underbrace{-\frac{1}{2}v\int(\partial_{x}^{l+1}u)^{2}(x, t)\chi_{n-l}^{'}(x+vt; \epsilon, b)\mbox{d}x}
_{\text{$A_{1}$}}\nonumber\\
&&+\frac{3}{2}\int(\partial_{x}^{l+2}u)^{2}(x, t)\chi_{n-l}^{'}(x+vt; \epsilon, b)\mbox{d}x
\underbrace{-\frac{1}{2}\int(\partial_{x}^{l+1}u)^{2}(x, t)\chi_{n-l}^{'''}(x+vt; \epsilon, b)\mbox{d}x}
_{\text{$A_{2}$}}\nonumber\\
&&\underbrace{+\int\partial_{x}^{l+1}(u\partial_{x}u)\partial_{x}^{l+1}u(x, t)\chi_{n-l}(x+vt; \epsilon, b)\mbox{d}x}
_{\text{$A_{3}$}}\nonumber\\
&&\underbrace{-\int H\partial_{x}^{l+3}u\partial_{x}^{l+1}u(x, t)\chi_{n-l}(x+vt; \epsilon, b)\mbox{d}x}
_{\text{$A_{4}$}}=0.
\end{eqnarray}
From the previous step $(n-l, l)$ , we derive


\begin{eqnarray*}\label{429}
\int_{\hat{t}}^{T}|A_{1}(t)|\mbox{d}t
\leq |v|\int_{\hat{t}}^{T}\int(\partial_{x}^{l+1}u)^{2}(x, t)\chi_{n-l}^{'}(x+vt; \epsilon, b)\mbox{d}x\mbox{d}t
\leq c_{3}.
\end{eqnarray*}

Invoking (\ref{prp6}), one obtains

\begin{eqnarray}\label{441}
|A_{2}(t)|
&\leq& c_{3}\int(\partial_{x}^{l+1}u)^{2}(x, t)\chi_{n-l-3}^{'}(x+vt; \epsilon, b)\mbox{d}x\nonumber\\
&+&c_{3}\int(\partial_{x}^{l+1}u)^{2}(x, t)\chi(x+vt; \epsilon/10, \epsilon/2)\mbox{d}x.
\end{eqnarray}
According to previous steps $(n-l-3, l+1)$ and $(0, l+1)$,
we know that (\ref{441}) is bounded.
Notice that the step $(0, l+1)$ is implied by the step
$(1, l)=(l+1-l, l)\leq (n-l, l)$.

For the term $A_{3}$, Leibniz formula leads to
\begin{eqnarray}
A_{3}
&=&d_{0}\int u\partial_{x}^{l+2}u\partial_{x}^{l+1}u\chi_{n-l}(x+vt; \epsilon, b)\mbox{d}x
+d_{1}\int \partial_{x}u(\partial_{x}^{l+1}u)^{2}\chi_{n-l}(x+vt; \epsilon, b)\mbox{d}x\nonumber
\end{eqnarray}
\begin{eqnarray}\label{5A3}
&+&d_{2}\int\partial_{x}^{2}u\partial_{x}^{l}u\partial_{x}^{l+1}u\chi_{n-l}(x+vt; \epsilon, b)\mbox{d}x
+\sum_{j=3}^{l-1}\int\partial_{x}^{j}u\partial_{x}^{l+2-j}u\partial_{x}^{l+1}u\chi_{n-l}(x+vt; \epsilon, b)\mbox{d}x\nonumber\\
&=&A_{3,0}+A_{3,1}+A_{3,2}+\sum_{j=3}^{l-1}A_{3, j}.
\end{eqnarray}
After integration by parts, we deduce
\begin{eqnarray*}
A_{3,0}
&=&-\frac{d_{0}}{2}\int\partial_{x}u(\partial_{x}^{l+1}u)^{2}\chi_{n-l}(x+vt; \epsilon, b)\mbox{d}x\nonumber\\
&-&\frac{d_{0}}{2}\int u(\partial_{x}^{l+1}u)^{2}\chi_{n-l}^{'}(x+vt; \epsilon, b)\mbox{d}x\nonumber\\
&=&A_{3,01}+A_{3,02}
\end{eqnarray*}
with $A_{3,01}$ be similar to $A_{3,1}$.

Sobolev embedding yields
\begin{eqnarray*}
\int_{\hat{t}}^{T}|A_{3, 02}(t)|\mbox{d}t
&\leq&\sup_{t\in[\hat{t}, T]}\|u\|_{L^{\infty}}
\int_{\hat{t}}^{T}\int(\partial_{x}^{l+1}u)^{2}\chi_{n-l}^{'}(x+vt; \epsilon, b)\mbox{d}x\mbox{d}t\nonumber\\
&\leq&\sup_{t\in[\hat{t}, T]}\|u\|_{H^{3/4^{+}}}
\int_{\hat{t}}^{T}\int(\partial_{x}^{l+1}u)^{2}\chi_{n-l}^{'}(x+vt; \epsilon, b)\mbox{d}x\mbox{d}t,
\end{eqnarray*}
where the last integral  corresponds to the case $(n-l, l)$, which
is part of our hypothesis of induction.

For the term $A_{3, 1}$, we have
\begin{eqnarray*}
|A_{3, 1}|
\leq
c_{3}\|\partial_{x}u\|_{L^{\infty}}\int(\partial_{x}^{l+1}u)^{2}\chi_{n-l}(x+vt; \epsilon, b)\mbox{d}x
\end{eqnarray*}
which can be handled by the Gronwall inequality.

Consider now $A_{3, 2}$ which appears only if $l\geq 2$ (we recall that $n\geq 3$(to be proved $(n-l, l+1))$)
\begin{eqnarray}\label{5A32}
A_{3, 2}
=d_{2}\int\partial_{x}^{2}u\partial_{x}^{l}u\partial_{x}^{l+1}u\chi_{n-l}(x+vt; \epsilon, b)\mbox{d}x.
\end{eqnarray}
Following the idea in \cite{ilp1},
we study two cases:
$l=2$ and $l\geq 3$.

We first consider $l=2$.

Similar to the estimates of (\ref{1A331})-(\ref{1A333}) in the previous section,
one derives
\begin{eqnarray}\label{448}
|A_{3, 2}|
&=&\left|-\frac{d_{2}}{3}\int\partial_{x}^{2}u\partial_{x}^{2}u\partial_{x}^{2}u\chi_{n-l}^{'}(x+vt; \epsilon, b)\mbox{d}x\right|\nonumber\\
&\leq&
c\|(\partial_{x}^{2}u)\chi_{n-l}^{'}(\cdot+vt; \epsilon, b)\|_{L^{\infty}}
\int(\partial_{x}^{2}u)^{2}\chi(x+vt; \epsilon/10, \epsilon/2)\mbox{d}x\nonumber\\
&\leq& c_{3}\|(\partial_{x}^{2}u)\chi_{n-l}^{'}(\cdot+vt; \epsilon, b)\|_{L^{\infty}}^{2}+c_{3}\nonumber\\
&\leq& c_{3}\|(\partial_{x}^{2}u)^{2}\chi_{n-l}^{'}(\cdot+vt; \epsilon, b)\|_{L^{\infty}}+c_{3}\nonumber\\
&\leq& c_{3}\int|\partial_{x}[(\partial_{x}^{2}u)^{2}\chi_{n-l}^{'}(x+vt; \epsilon, b)]|\mbox{d}x
+c_{3}\nonumber\\
&\leq& c_{3}\int|(\partial_{x}^{2}u)^{2}\chi_{n-l}^{'}(x+vt; \epsilon, b)|\mbox{d}x
+c_{3}\int|(\partial_{x}^{3}u)^{2}\chi_{n-l}^{'}(x+vt; \epsilon, b)|\mbox{d}x\nonumber\\
&+&c_{3}\int|(\partial_{x}^{2}u)^{2}\chi_{n-l}^{''}(x+vt; \epsilon, b)|\mbox{d}x\nonumber\\
&=&A_{3,21}+A_{3,22}+A_{3,23}+c_{3}.
\end{eqnarray}
From our induction hypothesis we know that $A_{3,21} $
and $A_{3,22}$ are bounded after integration in time,
since $A_{3,21} $  corresponds to
the case $(n-l, 1)=(n-2, 1)$ and
$A_{3,22} $ corresponds to the case $(n-l, 2)=(n-2, 2)$.

Moreover, invoking (\ref{prp5}), one derives
\begin{eqnarray}\label{5A323}
|A_{3,23}|
&\leq& c_{3}\int(\partial_{x}^{2}u)^{2}\chi_{n-l-2}(x+vt; \epsilon, b)\mbox{d}x
+c_{3}\int(\partial_{x}^{2}u)^{2}\chi^{'}(x+vt; \epsilon/3, b+\epsilon)\mbox{d}x\nonumber\\
&=&A_{3,231}+A_{3,232}.
\end{eqnarray}
From the induction cases $(n-l-2, 2)$ and $(0, 1)$,
we deduce that $A_{3,231}$  is bounded in
 time $t\in [\hat{t}, T]$  and
$A_{3,232}$ can be controlled after integration in time.

This completes the proof of (\ref{5A32}) in the case $l=2$.

Next, we turn to the case  $l\geq3.$

Integration by parts leads to
 \begin{eqnarray}\label{452}
 A_{3,2}
 &=&
d_{2}\int\partial_{x}^{2}u\partial_{x}^{l}u\partial_{x}^{l+1}u\chi_{n-l}(x+vt; \epsilon, b)\mbox{d}x\nonumber\\
&=&-\frac{d_{2}}{2}\int\partial_{x}^{3}u(\partial_{x}^{l}u)^{2}\chi_{n-l}(x+vt; \epsilon, b)\mbox{d}x\nonumber\\
&+&\frac{d_{2}}{2}\int\partial_{x}^{2}u(\partial_{x}^{l}u)^{2}\chi_{n-l}^{'}(x+vt; \epsilon, b)\mbox{d}x
\end{eqnarray}
 For the integrals on the right hand side
of (\ref{452}), using (\ref{433}) and reasoning as (\ref{448}) produce
\begin{eqnarray*}\label{w1}
 |A_{3,2}|
& \leq&
 c_{3}\int|(\partial_{x}^{2}u)^{2}\chi_{n-l}^{'}(x+vt; \epsilon, b)|\mbox{d}x
+c_{3}\int|(\partial_{x}^{3}u)^{2}\chi_{n-l}^{'}(x+vt; \epsilon, b)|\mbox{d}x\nonumber\\
&+&c_{3}\int|(\partial_{x}^{2}u)^{2}\chi_{n-l}^{''}(x+vt; \epsilon, b)|\mbox{d}x
+c_{3}\int|(\partial_{x}^{4}u)^{2}\chi_{n-l}^{'}(x+vt; \epsilon, b)|\mbox{d}x\nonumber\\
&+&c_{3}\int|(\partial_{x}^{3}u)^{2}\chi_{n-l}^{''}(x+vt; \epsilon, b)|\mbox{d}x.
\end{eqnarray*}
Since $l\geq 3,$ after integration in time,
the first two and the fourth integrals  correspond to the previous cases $(n-l, 1)$
, $(n-l, 2)$ and $(n-l, 3)$, respectively, which are all implied in the case $(n-l, l)$.
The third and fifth integrals  can be treated using a similar way as
(\ref{5A323}), where the fifth integral corresponds to the case
$(n-l-2, 3)$ and $(0, 2)$ after using (\ref{prp5}). Note
that the case $(0, 2)$ is implied by the case $(1, 1)$,
which can be deduced from   the previous case $l=1.$

Therefore, we only need to consider the remainder terms in (\ref{5A3}),
i.e.,
\begin{eqnarray*}\label{453}
 A_{3,j}
 =
c_{j}\int\partial_{x}^{j}u\partial_{x}^{l+2-j}u\partial_{x}^{l+1}u\chi_{n-l}(x+vt; \epsilon, b)\mbox{d}x.
\end{eqnarray*}

Without loss of generality , we can assume $3\leq j \leq l/2+1.$
 Consequently, one finds
 \begin{eqnarray*}\label{454}
 |A_{3, j}|
\leq
c_{j}\int(\partial_{x}^{j}u\partial_{x}^{l+2-j}u)^{2}\chi_{n-l}(x+vt; \epsilon, b)\mbox{d}x
+c_{j}\int(\partial_{x}^{l+1}u)^{2}\chi_{n-l}(x+vt; \epsilon, b)\mbox{d}x
\end{eqnarray*}
with the second integral be the quantity to be estimated.

For the first integral, we have
\begin{eqnarray}\label{455}
&&c_{j}\int(\partial_{x}^{j}u\partial_{x}^{l+2-j}u)^{2}\chi_{n-l}(x+vt; \epsilon, b)\mbox{d}x\nonumber\\
&&\leq \|(\partial_{x}^{j}u)^{2}\chi(x+vt; \epsilon/10, \epsilon/2)\|_{L^{\infty}}
\int(\partial_{x}^{l+2-j}u)^{2}\chi_{n-l}(x+vt; \epsilon, b)\mbox{d}x.
\end{eqnarray}
From the induction hypothesis $(n-l, l+2-j)$ with $j\geq3$, we deduce that
the integral in (\ref{455}) is bounded.
Thus it remains to control the $L^{\infty}$ norm.

For this purpose, we employ the Sobolev inequality
$\|f\|_{L^{\infty}}\leq \|f\|_{H^{1, 1}}$ to obtain
\begin{eqnarray*}\label{456}
&&\|(\partial_{x}^{j}u)^{2}\chi(x+vt; \epsilon/10, \epsilon/2)\|_{L^{\infty}}\nonumber\\
&&\leq \int|\partial_{x}[(\partial_{x}^{j}u)^{2}\chi(x+vt; \epsilon/10, \epsilon/2)]|\mbox{d}x\nonumber\\
&&\leq c\int|\partial_{x}^{j}u\partial_{x}^{j+1}u\chi(x+vt; \epsilon/10, \epsilon/2)|\mbox{d}x
+c\int|(\partial_{x}^{j}u)^{2}\chi^{'}(x+vt; \epsilon/10, \epsilon/2)|\mbox{d}x\nonumber\\
&&\leq c\int(\partial_{x}^{j}u)^{2}\chi(x+vt; \epsilon/10, \epsilon/2)\mbox{d}x
+c\int(\partial_{x}^{j+1}u)^{2}\chi(x+vt; \epsilon/10, \epsilon/2)\mbox{d}x\nonumber\\
&&+c\int|(\partial_{x}^{j}u)^{2}\chi^{'}(x+vt; \epsilon/10, \epsilon/2)|\mbox{d}x.
\end{eqnarray*}

Since $j\leq l-1$, we have $j+1\leq l\leq n.$
Thus, previous cases $(0, j)$ and
$(0, j+1)$ imply the boundedness of the first two integrals, respectively.
The third integral corresponds to the case $(0, j-1)$ after integration in time.

Finally, we estimate $A_{4}.$
As before, we write
\begin{eqnarray*}
A_{4}
&=&-\int H\partial_{x}^{l+3}u\partial_{x}^{l+1}u\chi_{n-l}(x+vt;  \epsilon, b)\mbox{d}x\nonumber\\
&=&\int H\partial_{x}^{l+2}u\partial_{x}^{l+2}u\chi_{n-l}(x+vt;  \epsilon, b)\mbox{d}x
+\int H\partial_{x}^{l+2}u\partial_{x}^{l+1}u\chi_{n-l}^{'}(x+vt;  \epsilon, b)\mbox{d}x\nonumber\\
&=&A_{41}+A_{42}.
\end{eqnarray*}
The term $A_{41}$ can be treated easily, we omit it.

For the term $A_{42},$ one has
\begin{eqnarray*}
A_{42}
&=&\int H\partial_{x}^{l+2}u\partial_{x}^{l+1}u\chi_{n-l}^{'}(x+vt)\mbox{d}x\nonumber\\
&=&\int H\partial_{x}^{l+2}u \eta_{n-l} \partial_{x}^{l+1}u\eta_{n-l}\mbox{d}x\nonumber\\
&=&\int  H(\partial_{x}^{l+2}u \eta_{n-l}) \partial_{x}^{l+1}u\eta_{n-l}\mbox{d}x
-\int [H; \eta_{n-l}]\partial_{x}^{l+2}u\partial_{x}^{l+1}u\eta_{n-l}\mbox{d}x\nonumber\\
&=&\int  H\partial_{x}(\partial_{x}^{l+1}u \eta_{n-l}) \partial_{x}^{l+1}u\eta_{n-l}\mbox{d}x
-\int  H(\partial_{x}^{l+1}u \eta_{n-l}^{'}) \partial_{x}^{l+1}u\eta_{n-l}\mbox{d}x\nonumber\\
&-&\int [H; \eta_{n-l}]\partial_{x}^{l+2}u\partial_{x}^{l+1}u\eta_{n-l}\mbox{d}x\nonumber\\
&=&A_{421}+A_{422}+A_{423},
\end{eqnarray*}
where $A_{421}$ is positive and
 will stay at the left hand side of (\ref{438}).

(\ref{prp10}) and (\ref{prp11}) lead to
\begin{eqnarray}\label{6A422}
|A_{422}|
&\leq&\left|\int  H(\partial_{x}^{l+1}u \eta_{n-l}^{'}) \partial_{x}^{l+1}u\eta_{n-l}\mbox{d}x\right|\nonumber\\
&\leq& c\int(\partial_{x}^{l+1}u)^{2} (\eta_{n-l}^{'})^{2}\mbox{d}x
+c\int(\partial_{x}^{l+1}u \eta_{n-l})^{2}\mbox{d}x\nonumber\\
&\leq& c\int(\partial_{x}^{l+1}u)^{2} \chi_{n-l-1}(x; \epsilon/3, b+\epsilon)\mbox{d}x,
\end{eqnarray}

which can be handled by the previous step $(n-l-1, l+1)$
since $l+1\leq n.$

The term $A_{423}$ can be handled similarly, we omit it.
%
%
%



This basically completes the proof of Theorem \ref{persistence}.

To justify the previous formal computation, we
approximate
the initial data $u_{0} $ by    Schwartz  functions,
say $u_{0}^{\mu}$, $\mu> 0,$  which can be  satisfied
by convolution  $u_{0} $ with a family of mollifiers.
Using the well-posedness in the  class of Schwartz  functions,
 we obtain a family of solutions $u^{\mu}(\cdot, t)$
for which each step of the above argument can be justified.
From our construction
those estimates are uniform in the parameter $\mu> 0$, which yields the desired
estimate by passing to the limit.

\textbf{Acknowledgement}

The research of B. Guo
is partially supported by the National Natural Science Foundation
of China, grant 11731014.



\begin{thebibliography}{}



\bibitem{abr}
   J.P. Albert, J.L. Bona, J.M. Restrepo,
   Solitary wave solutions of the Benjamin equation,
    \textit{SIAM J. Appl. Math.}
    59 (1999) 2139-2161.



\bibitem{cb}
    H. Chen, J.L. Bona,
    Existence and asymptotic properties of solitary-wave solutions of Benjamin-type equations,
    \textit{Adv.  Differential Equations}
    3 (1998) 51-84.


\bibitem{be}
   T.B. Benjamin,
    A new kind of solitary wave,
    \textit{J. Fluid Mech.}
    245 (1992) 401-411.




\bibitem{bo}
   J. Bourgain,
    Fourier transform restriction phenomena for certain lattice subsets and applications to nonlinear evolution equations,
    \textit{Geom. Funct. Anal. }
     3 (1993) 107-156.




\bibitem{ca}
    A.P. Calder\'{o}n,
    Commutators of singular integral operators,
    \textit{Proc. Natl. Acad. Sci.}
    340 (1965) 1092-1099.










\bibitem{cgx}
   W. Chen, Z. Guo, J. Xiao,
    Sharp well-posedness for the Benjamin equation,
    \textit{Nonlinear Anal.}
      74 (2011) 6209-6230.


\bibitem{dmp}
    L. Dawson, H. McGahagan, G. Ponce,
    On the decay properties of solutions to a class of Schr\"{o}dinger equations,
    \textit{Proc. Amer. Math. Soc.}
     136 (2008) 2081-2090.

\bibitem{fp1}
    G. Fonseca, G. Ponce,
    The IVP for the Benjamin-Ono equation in weighted Sobolev spaces,
    \textit{J. Funct. Anal.}
     260 (2011) 436-459.


\bibitem{flp1}
    G. Fonseca, F. Linares, G. Ponce
    The IVP for the Benjamin-Ono equation in weighted Sobolev spaces II,
    \textit{J. Funct. Anal.}
      262 (2012) 2031-2049.



\bibitem{flp2}
    G. Fonseca, F. Linares, G. Ponce
    The IVP for the dispersion generalized Benjamin-Ono equation in weighted Sobolev spaces,
    \textit{Ann. Inst. H. Poincar\'{e} Anal. Non Lin\'{e}aire}
      30 (2013) 763-790.





\bibitem{gh2}
    B. Guo, Z. Huo,
    The cauchy problem for the generalized Korteweg-de Vries-Benjamin-Ono equation with low regularity data,
    \textit{Acta Math. Sin. (Engl. Ser.)}
      21 (2005) 1191-1196.




\bibitem{ilp1}
    P. Isaza, F. Linares, G. Ponce,
    On the propagation of regularity and decay of solutions to the k-generalized Korteweg-de Vries equation,
    \textit{Comm. Partial Differential Equations}
    40  (2015) 1336-1364.


\bibitem{ilp2}
    P. Isaza, F. Linares, G. Ponce,
    On the propagation of regularities in solutions of the Benjamin-Ono equation,
    \textit{J. Funct. Anal.}
    270 (2016) 976-1000.








\bibitem{ka}
    T. Kato,
    On the Cauchy problem for the (generalized) Korteweg-de Vries equation, in:
    \textit{ Adv.
      Math. Suppl. Stud., Stud. Appl. Math.}
    8 (1983) 93-128.

\bibitem{kp}
    T. Kato and G. Ponce,
    Commutator estimates and the Euler and Navier-Stokes equations,
    \textit{Comm.  Pure  Appl. Math.}
    41  (1988) 891-907.


\bibitem{kpv}
    C. Kenig, G. Ponce, L. Vega,
    A bilinear estimate with applications to the KdV equation,
    \textit{J. Amer. Math. Soc.}
    9 (1996) 573-603.




\bibitem{la}
    C. Laurey,
    The cauchy problem for a third order nonlinear Schr\"{o}dinger equation,
    \textit{Nonlinear Anal.}
     29  (1997) 121-158.

\bibitem{li}
    F. Linares,
    $L^{2}$ Global well-posedness of the initial value problem associated to the Benjamin equation,
    \textit{J.  Differential Equations}
     152 (1999) 377-393.

 \bibitem{lw}
    Y. Li, Y. Wu,
    Global well-posedness for the Benjamin equation in low regularity,
    \textit{Nonlinear Anal.}
       73  (2010) 1610-1625.

\bibitem{np}
    J. Nahas, G. Ponce,
    On the persistent properties of solutions to semi-linear Schr\"{o}dinger equation,
    \textit{Comm. Partial Differential Equations}
     34 (2009) 1-20.

\bibitem{pa}
   J.A. Pava,
  Existence and stability of solitary wave solutions of the Benjamin equation,
    \textit{J.  Differential Equations}
   152  (1999) 136-159.

\bibitem{sl}
    M. Scialom, F. Linares,
    On generalized Benjamin type equations,
    \textit{Discrete Contin. Dyn. Syst.}
    12  (2005) 161-174.


\bibitem{ss}
    J.I. Segata, D.L. Smith,
    Propagation of regularity and persistence of decay for fifth order dispersive models,
    \textit{J. Dynam. Differential Equations}
    (2015) 1-36.

\bibitem{u}
    J.J. Urrea,
     The Cauchy problem associated to the Benjamin equation in weighted Sobolev spaces,
    \textit{J.  Differential Equations}
     254  (2013) 1863-1892.













\end{thebibliography}
\end{document}